\documentclass[final,leqno]{siamltex}

\usepackage{amsmath,amssymb}
\usepackage{mathtools}
\usepackage{hyperref} 
\hypersetup{
 colorlinks=true,
 linktocpage=true,
 bookmarksnumbered=true
}
\usepackage{algorithm,algorithmic}
\usepackage[pdftex]{graphicx,color}

\newtheorem{thm}{Theorem}[section]
\newtheorem{rem}[thm]{Remark}

\newcommand{\seqof}[3]{(#1)_{#2}^{#3}}

\newcommand{\ssa}{\text{SSA}}

\newcommand{\MNRM}{\text{MNRM}}
\newcommand{\TL}{\text{TL}}

\newcommand{\NMNRM}{N_{\MNRM}}
\newcommand{\NMNRMC}{N_{\MNRM}^{(c)}}

\newcommand{\NSSAP}{N_{{\ssa}^*}}
\newcommand{\NTL}{N_{\TL}}

\newcommand{\NTLC}{N_{\TL}^{(c)}}

\newcommand{\EE}{\mathcal{E}}

\newcommand{\WEH}[1]{\hat \EE_{I,#1}}

\newcommand{\avg}[2]{\mathcal{A}\left(#1;#2\right)}
\newcommand{\svar}[2]{\mathcal{S}^2\left(#1;#2\right)}

\newcommand{\dbar}[1]{\Bar{\Bar{#1}}}

\newcommand{\ie}{\emph{i.e.}}
\newcommand{\eg}{\emph{e.g.}}

\newcommand{\ud}{\mathrm{d}}
\newcommand{\prob}[1]{\mathrm{P}\left(#1\right)}
\newcommand{\expt}[1]{\mathrm{E}\left[#1\right]}

\newcommand{\var}[1]{\mathrm{Var}\left[#1\right]}

\newcommand{\rset}{\mathbb{R}}

\newcommand{\zset}{\mathbb{Z}}
\newcommand{\JAC}{\mathbb{J}}

\newcommand{\Ordo}[1]{{\mathcal{O}}\left(#1\right)}

\newcommand{\ordo}[1]{{o}\left(#1\right)}

\newcommand{\latt}{\mbox{$\zset_+^d$}}

\newcommand{\indicator}[1]{\mathbf{1}_{#1}} 

\newcommand{\PERIOD}{.}
\newcommand{\COMMA}{,}

\newcommand{\LP}{\left(}
\newcommand{\RP}{\right)}

\newcommand{\SEP}{\, \big| \,}
\newcommand{\BX}{\bar X(t)}
\newcommand{\BXi}{\bar X_i(t)}

\newcommand{\sj}{\sum_{j=1}^J}

\newcommand{\mV}[1]{\mathcal{V}_{#1}}
\newcommand{\hmV}[1]{\hat{\mathcal{V}}_{#1}}

\title{A multilevel adaptive reaction-splitting simulation method for stochastic reaction networks}

\author{Alvaro Moraes\thanks{Mathematical and Computer Sciences and Engineering Division,
 King Abdullah University of Science and Technology (KAUST),
 Thuwal, Saudi Arabia ({\tt alvaro.moraesgutierrez@kaust.edu.sa}).}
        \and Raul Tempone\thanks{Mathematical and Computer Sciences and Engineering Division,
 King Abdullah University of Science and Technology (KAUST),
 Thuwal, Saudi Arabia ({\tt raul.tempone@kaust.edu.sa}).}
\and Pedro Vilanova\thanks{Mathematical and Computer Sciences and Engineering Division,
 King Abdullah University of Science and Technology (KAUST),
 Thuwal, Saudi Arabia ({\tt pedro.guerra@kaust.edu.sa}).}}
 
\numberwithin{equation}{section}
\numberwithin{figure}{section}
\numberwithin{table}{section}

\begin{document}
\maketitle
\begin{abstract}
Stochastic modeling of reaction networks is a framework used to describe the time evolution of many natural and artificial systems, including, biochemical reactive systems at the molecular level, viral kinetics, the spread of epidemic diseases, and wireless communication networks, among many other examples.
In this work, we present a novel multilevel Monte Carlo method for kinetic simulation of stochastic reaction networks  that is specifically designed for systems in which the set of reaction channels can be adaptively partitioned into two subsets characterized by either ``high'' or ``low'' activity. 
Adaptive in this context means that the partition evolves in time according to the states visited by the stochastic paths of the system.
To estimate expected values of observables of the system at a prescribed final time, our method bounds the global computational error to be below a prescribed tolerance, $TOL$, within a given confidence level. This is achieved with a  computational complexity of order $\Ordo{TOL^{-2}}$, the same as with an exact method, but with a smaller constant. We also present a novel control variate technique based on the stochastic time change representation by Kurtz, which may dramatically reduce the variance of the coarsest level
at a negligible computational cost.
Our numerical examples show substantial gains with respect to the standard Stochastic Simulation Algorithm (SSA) by Gillespie and also our previous hybrid Chernoff tau-leap method. 
\end{abstract}



\setcounter{tocdepth}{1}

\section{Introduction}
Stochastic reaction networks (SRN) are mathematical models that employ Markovian dynamics to describe the {time evolution of interacting particle systems} where one particle interact with the others through a finite set of reaction channels.
Typically, there is a finite number of interacting chemical species $(S_1,S_2,\ldots,S_d)$ and a stochastic process, $X$, such that its $i$-th coordinate is a non-negative integer number $X_i(t)$ that keeps track of the abundance of the $i$-th species at time $t$.  Therefore, the state space of the process $X$ is the lattice $\latt$.

Our main goal is to estimate the expected value $\expt{g(X(T))}$, where $X$ is a non-homogeneous Poisson process describing a SRN, and $g:\rset^d\to\rset$ is a given real observable of $X$ at a final time $T$. 
Pathwise realizations can be simulated
exactly using the {Stochastic Simulation Algorithm} (SSA),
introduced by Gillespie in \cite{gillespie_ssa} (also known as Kinetic Monte Carlo among physicists, see \cite{KMC} and references therein), or the {Modified Next Reaction Method} (MNRM) introduced by Anderson in \cite{anderson2007modified}, among other methods. Although these algorithms generate exact realizations of $X$, they may be  computationally expensive for systems that undergo high activity. For that reason, Gillespie proposed in
\cite{gillespie_tau_leap} the tau-leap method to approximate the SSA by evolving the process with fixed time
steps while freezing the propensity functions at the beginning of each time step. 

A drawback of the tau-leap method is that the simulated paths may take negative values, which is a nonphysical consequence of the approximation and not a qualitative feature of the original process. For that reason, we proposed in \cite{ourSL,ourML} a Chernoff-based hybrid method that switches adaptively between the tau-leap and an exact method. This allows us to control the probability of reaching negative values while keeping the computational work substantially smaller than the work of an exact method.
The hybrid method developed in \cite{ourSL,ourML} can be successfully applied  to systems where the state space, $\zset_+^d$, can be decomposed into two regions   according to the activity of the system; where all the propensities are uniformly low or uniformly high, \ie, non-stiff systems.
To handle stiff systems, we first measure the total activity of the system at a certain state by the total sum of the propensity functions evaluated at this state. 
The activity of the system is low when all the propensities are uniformly low, but a high level of activity can be the result of a high activity level in one single channel.   
This observation suggests that to reduce computational costs, we should adaptively split the set of reaction channels into two subsets according to the individual high and low activity levels.  
It is natural to evolve the system in time by applying the tau-leap method to the high activity channels and an exact method to the low activity ones. 
This is the main idea we develop in this work.

Reaction-splitting methods for simulating stochastic reaction networks are treated for instance in \cite{Harris2006, Puchalka2004, Haseltine2002, Plyasunov2005}, but our work is, to the best of our knowledge, the first that  
i) achieves the computational complexity of an exact method like the SSA by using the multilevel Monte Carlo paradigm, 
ii) explicitly uses a  decomposition of the global error to provide all the simulation parameters needed to achieve our goal with minimal computational effort,
iii) effectively controls the global probability of reaching negative populations with the tau-leap method, and 
iv) needs only two user-defined parameters that are natural quantities - the maximum allowed relative global error or tolerance and the confidence level.  

In \cite{Harris2006}, the authors propose an adaptive reaction-splitting scheme that considers not only the exact and tau-leap methods but also the Langevin and Mean Field ones. 
Their main goal is to obtain fast hybrid simulated paths, and they do not try to control the global error. 
The efficiency of their method is measured a posteriori using smoothed frequency histograms that should be close to the exact ones according to the distance defined by Cao and Petzold in \cite{CaoPetzold2006}. 
In their work, the tau-leap step is chosen according to the ``leap condition'' (as in \cite{Cao2006}) but they do not perform a rigorous control of the global discretization error.
In order to avoid negative populations, the authors reverse population updates if any value is found to be negative after accounting for all the reactions. Then, the tau-lep step size is decremented and the path simulation is restarted. 
This approach introduces bias in the estimations, and even by controlling the small reactant populations, a tau-leap step always may lead to negative populations subsequently increasing its computational work. Our Chernoff -based bound is a fast and accurate procedure to obtain the correct tau-leap step size. 
Finally, the method in \cite{Harris2006} needs to define three parameters that quantify the speed of the reaction channels, which, in principle, are not trivial to determine for a given problem.

Puchalka and Kierzek's approach \cite{Puchalka2004} seems to be closest to our approach in spirit since they also explore the idea of adaptively splitting the set of reaction channels using the tau-leap method for the fast ones and an exact method for the slow ones. 
They seek to simulate fast approximate paths while maintaining qualitative features of the system. 
The quantitative features are checked a posteriori against an exact method.  
Regarding their tau-leap step size selection, Puchalka and Kierzek consider a user-defined maximal time step empirically chosen by numerical tests instead of controlling the discretization error.
Their classification rule is applied individually to each reaction channel.  It takes into account both the percentage of individual activity and the abundance of the species consumed. In a certain sense it can be seen as a way of controlling the probability of negative populations and an ad-hoc manner to split the reaction channels by optimizing the computational work.

In \cite{Haseltine2002} and \cite{Plyasunov2005}, the reaction-splitting issue is addressed but the partition method is not adaptive, \ie, fast and slow reaction channels are identified offline and are inputs of the algorithms.
We note that these works do not provide any measure or control of the resulting global error. 
 Furthermore, they do not control the probability of attaining negative populations.   

In the remaining of this section, we introduce the mathematical model and the path simulation techniques used in this work.
In Section \ref{sec:mixedpaths}, we present an algorithm to generate mixed trajectories; that is, the algorithm generates a trajectory using an exact method for the low activity channels and the Chernoff tau-leap method for the high activity ones. Then, inspired by the ideas of Anderson and Higham \cite{Anderson2012}, we propose an algorithm for coupling two mixed Chernoff tau-leap paths. This algorithm uses four building blocks that result from the combination of the MNRM and the tau-leap methods.
In Section \ref{sec:mlmcerr}, we propose a mixed MLMC estimator. Next, we introduce a global error decomposition and show that the computational complexity of our method is of order $\Ordo{TOL^{-2}}$.
Finally, we show the automatic procedure that estimates our quantity of interest within a given prescribed relative tolerance, up to a given confidence level.
Next, in Section \ref{sec:cvar}, we present a novel control variate technique to reduce the variance of the quantity of interest at level 0. 
In Section \ref{sec:examples}, the numerical examples illustrate the advantages of the mixed MLMC method over the hybrid MLMC method presented in \cite{ourML} and to the SSA. Finally, Section \ref{sec:conclusions} presents our conclusions.

\subsection{A Class of Markovian Pure Jump Processes}
\label{sec:pjp}
In this section, we describe the class of Markovian pure jump processes, $X:[0,T]\times \Omega \to\zset_+^d$, frequently used for modeling stochastic biochemical reaction networks.  
 
Consider a biochemical system of $d$ species interacting
through $J$ different reaction channels. For the sake of brevity, we write $X(t,\omega) {\equiv} X(t)$.  Let $X_i(t)$ be the
number of particles of species $i$ in the system at time $t$. We 
study the evolution of the state vector, $X(t) = (X_1(t), \ldots, X_d(t)) \in
  \zset_+^d$,
modeled as a continuous-time Markov chain starting at $X(0) \in \zset_+^d$. Each reaction can be described by the
vector $\nu_j\in \zset^d$, such that, for a state vector $x\in
\zset_+^d$, a single firing of reaction $j$ leads to the change
$x \to x+ \nu_j$.
The probability that reaction $j$ will occur during the small interval
$(t,t{+}\ud t)$ is then assumed to be
\begin{equation}\label{eq:prob}
  \prob{X(t+\ud t)=x+\nu_j | X(t) = x} =  a_j(x) \ud t + \ordo{\ud t}
\end{equation}
for a given non-negative {polynomial propensity function}, $a_j:\rset^d \to
\rset$. We set $a_j(x){=}0$ for those $x$ such that $x{+}\nu_j\notin\latt$.
The process $X$ admits the following random time change representation by Kurtz \cite{kurtzmp}:
\begin{equation}
  \label{eq:exact_process}
  X(t) = X(0) + \sum_{j=1}^J  \nu_j Y_j \LP \int_0^t a_j(X(s))\, \ud s \RP  \COMMA
\end{equation}
where $Y_j:\rset_+ {\times} \Omega \to\zset_+$ are independent unit-rate Poisson processes.
Hence, $X$ is a non-homogeneous Poisson process.

\subsection{The Modified Next Reaction Method (MNRM)}\label{sec:MNR}
The MNRM, introduced in \cite{anderson2007modified}, and based on the Next
Reaction Method  \cite{GibBruck}, is an exact simulation algorithm like Gillespie's SSA that explicitly uses  representation \eqref{eq:exact_process} for simulating exact paths and generates only one exponential random variable per iteration. The reaction times are modeled with firing times of Poisson processes, $Y_j$, with
internal times given by the integrated propensity functions.
The randomness is now
separated from the state of the system and is encapsulated in the $Y_j$'s.
Computing the next reaction and its time is equivalent to computing how
much time passes before one of the Poisson process, $Y_j$, fires, and
which process fires at that particular time, by taking the minimum of such times. 
It is important to mention that the MNRM is used to simulate correlated exact/tau-leap paths as well as nested  tau-leap/tau-leap paths, as in \cite{ourML,Anderson2012}.
In Section \ref{cmp}, we use this feature for coupling two mixed paths.

\subsection{The Tau-Leap Approximation}\label{sec:num_approx}

In this section, we define $\bar X$, the tau-leap approximation of the process, $X$, which follows from applying the forward-Euler approximation to the integral term in the random time change representation \eqref{eq:exact_process}.

The tau-leap method was proposed
in \cite{gillespie_tau_leap} to avoid the computational drawback of the exact methods, \ie,  when many reactions occur
during a short time interval. The tau-leap process, $\bar X$, starts from $X(0)$ at time $0$, and 
given that $\bar X(t) {=} \bar x$ and a
time step $\tau{>}0$, we have that
$\bar X$ at time $t{+}\tau$ is generated by
\begin{equation*}
  \bar X(t+\tau) = \bar x + \sum_{j=1}^J 
  \nu_j \mathcal{P}_j\left( a_j(\bar x)\tau \right) \COMMA
\end{equation*}
where $\{\mathcal{P}_j(\lambda_j)\}_{j=1}^J$ are independent
Poisson distributed random variables with parameter $\lambda_j$,
used to model the number of times that the reaction $j$ fires
during the $(t,t{+}\tau)$ interval.  Again, this is nothing else than a 
forward-Euler discretization of the stochastic differential equation formulation of the pure jump process
\eqref{eq:exact_process}, realized by the Poisson random measure with
state-dependent intensity (see, \eg, \cite{li}).

In the limit, when $\tau$ tends to zero, the tau-leap method gives the same solution
as the exact methods \cite{li}. The total number of firings in each channel is 
a Poisson distributed stochastic variable depending only on the
initial population, $\bar X(t)$.  The error thus comes from the
variation of $a(X(s))$ for $s\in(t,t{+}\tau)$. 

\subsection{The Hybrid Chernoff Tau-leap Method}
\label{sec:chernoff}
In \cite{ourSL}, we derived a Chernoff-type bound that allows us to guarantee that
the one-step exit probability in the tau-leap method is less than a
predefined quantity, $\delta{>}0$. 
The idea is to find the largest possible time step, 
$\tau$, such that, with high probability, in the next step, the approximate process, $\bar{X}$, will take a value in the lattice, $\latt$, of non-negative integers.
%
This can be achieved by solving $d$ auxiliary
problems, one for each $x$-coordinate, $\bar X_i(t)$,
$i=1,2,\ldots,d$ as follows. Find the largest possible $\tau_i\geq 0$, such that
\begin{equation}
\label{eq:taui-deltai}
\prob{  \BXi + \sj \nu_{ji} \mathcal{P}_j\LP a_j\LP \bar X(t)\RP\tau_i\RP < 0 \Biggm| \BX}\leq \delta_i \COMMA
\end{equation}
where $\delta_i {=} \delta/d$, and $\nu_{ji}$ is the $i$-th coordinate of
the $j$-th reaction channel, $\nu_j$. Finally, we let $\tau{:=}\min\{\tau_i: i=1,2,\ldots,d\} \PERIOD$
Using the exact pre-leap method we developed in \cite{ourSL,ourML} for single-level and multilevel hybrid schemes, allows us to switch adaptively between the tau-leap and an exact method.
By construction, the probability that one hybrid path exits the lattice, $\latt$, can be estimated by
\begin{align*}
\prob{ A^c} \leq \expt{1-(1-\delta)^{N_{\TL}}} = \delta \expt{N_{\TL}} - \frac{\delta^2}{2}(\expt{N_{\TL}^2} - \expt{N_{\TL}}) + o(\delta^2),
\end{align*}
where  $\bar{\omega} \in A$ if and only if the whole hybrid path, $\seqof{\bar X(t_k,\bar \omega)}{k=0}{K(\bar{\omega})}$, belongs to the lattice, $\latt$, $\delta {>}0$ is the one-step exit probability bound, and $N_{\TL}(\bar \omega) {\equiv} N_{\TL}$ is the number of tau-leap steps in a hybrid path. Here, $A^c$ is the complement of the set $A$.

To simulate a hybrid path, given the current state of the approximate process, $\BX$, we adaptively determine whether to use an exact or the tau-leap method  for the next step. This decision is based on the relative computational cost of taking an exact step versus the cost of taking a Chernoff tau-leap step.
Instead, in the present work, at each time step, we adaptively determine which reactions are suitable for using the exact method and which reactions are suitable for the Chernoff tau-leap method.


\section{Generating Mixed Paths}
\label{sec:mixedpaths}
In this section we explain how mixed paths are generated. First, we present the splitting heuristic; that is, we discuss how to partition the set of reaction channels at each decision time. Then, we present the one-step mixing rule, which is the main building block for constructing a mixed path. Finally, we show how to couple two mixed paths.

\subsection{The Splitting Heuristic}
\label{sec:splitting}
In this section, we explain how we partition the set of reaction channels, $\mathcal{R}{:=}\{1,...,J\}$, into $\mathcal{R}_{\TL}$ and $\mathcal{R}_{\MNRM}$.

Let $(t,x)$ be the current time and state of the approximate process, $\bar X$, and $H$  be the next decision (or synchronization) time (given by the Chernoff tau-leap step size $\tau_{Ch}=\tau_{Ch}(x,\delta)$ and the time mesh).
We want to split $\mathcal{R}$ into two subsets, $\mathcal{R}_{\MNRM}$ and $\mathcal{R}_{\TL}$, such that 
the expected computational work of reaching $H$, starting at $t$, is minimal for all possible splittings.

The idea goes as follows. First, we define a  linear order on $\mathcal{R}$, based on the basic principle that we want to use tau-leap for the $j$-th reaction if its activity is high. This linear order determines $J{+}1$ possible splittings, out of $2^J$. 
In order to measure the activity, it turns out that using only the propensity functions evaluated at $x$, that is, $a_j(x)$, is not enough. This is because the $j$-th reaction could affect components of $x$ with small values. If this is the case, this determines small Chernoff tau-leap step sizes. In order to avoid this scenario, we penalize the $j$-th reaction channel if it has a high exit probability. We approximate this exit probability using a Poisson distribution for each dimension of $x$. For example, let $\nu_{ji}$ be the $i$-th component of the $j$-th reaction channel. If $\nu_{ji}<0$, then the probability that a Poisson distributed random variable with rate $a_j(x)(H{-}t)$ is greater than $x_i/\nu_{ji}$ measures how likely species $x_i$ can become negative in the interval $H{-}t$, independently of reactions $j'{\in} \mathcal{R}$, $j{\neq} j'$.
Let $I_j{:=} \{i:\nu_{ji} <0\}$,
\begin{align}\label{eq:thetas}
\theta_{j}:= \left\{ 
\begin{array}{ll}
\prob{\mathcal{P}(a_j(x)(H{-}t)) > \min_{i\in I_j} \{-\frac{x_i}{\nu_{ji}}\}\SEP x} & \text{if } I_j \neq \emptyset\\
0 &\text{otherwise}
\end{array}
\right. .
\end{align}
Then, the penalty weight for $a_j(x)$ is $1{-}\theta_j$. We define $\tilde{a}_j(x){:=}(1{-}\theta_j)a_j(x)$. The linear order is then a permutation, $\sigma$, over $\mathcal{R}$ such that
$$\tilde{a}_{\sigma(j)}(x)>\tilde{a}_{\sigma(j+1)}(x), \,\,\,  j{=}1,...,J{-}1 \PERIOD$$
 
Second, we find among the $J{+}1$ partitions the one with optimal work. This  is the computational work incurred when performing one step of the algorithm using tau-leap for the reactions $\mathcal{R}_{TL}$ and the MNRM for the reactions $\mathcal{R}_\MNRM$. The work corresponding to $\mathcal{R}_{\TL}$ is 
\begin{align}\label{eq:worktl}
\text{Work}(\mathcal{R}_{\TL},x,t):=\frac{H{-}t}{\min\{\tau_{Ch},H{-}t\}} \left(C_{s} + \sum_{j\in \mathcal{R}_{\TL}}C_P(a_{j}(x) \tau_{Ch})\right) \COMMA
\end{align}
where $C_s$ is the work of computing the split (see Section \ref{sec:cs}), and $C_P(\lambda)$ is the work of a Poisson random variate with rate $\lambda$. 
The factor $\frac{H{-}t}{\min\{\tau_{Ch},H{-}t\}}$ takes into account the number of steps required to reach $H=H(t)$ from $t$.
For the Gamma simulation method developed by Ahrens and Dieter in \cite{poivariate}, which is the one used by MATLAB, $C_P$ is defined as
\begin{align*}
C_P(\lambda):=
\left\{ \!
\begin{array}{ll}
b_1{+}b_2 \ln \lambda & \text{ for } \lambda >15 \\
b_3{+}b_4 \lambda & \text{ for } \lambda \leq 15
\end{array}
\right. \PERIOD
\end{align*}
In practice, it is possible to estimate $b_i$, $i{=}1,2,3,4$ using Monte Carlo sampling and a least squares fit. 
For more details, we refer to \cite{ourSL}.

Similarly, the work corresponding to $\mathcal{R}_\MNRM$ is
$$\text{Work}(\mathcal{R}_{\MNRM},x,t):=\frac{H{-}t}{\min\{\tau_{\MNRM},H{-}t\}} \, C_{\MNRM} \COMMA$$
where the constant $C_{\MNRM}$ is the work of an MNRM step and $\tau_{\MNRM}{=}\LP \sum_{j\in \mathcal{R}_\MNRM} a_j(x)\RP^{-1}$.
\subsection{On the Work required to the Splitting Heuristic, $C_s{=}C_s(J)$}
\label{sec:cs}
The work required to perform the splitting includes the work required to determine $\text{Work}(\mathcal{R}_{\TL})$ and $\text{Work}(\mathcal{R}_{\MNRM})$, both defined in Section \ref{sec:splitting}. The linear order previously defined determines $J{+}1$ possible splittings, $\mathcal{S}_i$, $i{=}0,...,J$, as follows:
\begin{center}
\begin{tabular}{lc|c}
& $\mathcal{R}_{\TL}$ & $\mathcal{R}_\MNRM$ \\
\hline 
$\mathcal{S}_0$ & $\emptyset$ & $\mathcal{R}$ \\
$\mathcal{S}_1$ & $\{\sigma^{-1}(1)\}$ & $\{\sigma^{-1}(2),..,\sigma^{-1}(J)\}$ \\
$\mathcal{S}_2$ & $\{\sigma^{-1}(1),\sigma^{-1}(2)\}$ & $\{\sigma^{-1}(3),..,\sigma^{-1}(J)\}$ \\
$\vdots$ & & \\
$\mathcal{S}_{J}$ & $\mathcal{R}$ & $\emptyset$
\end{tabular} .
\end{center}
The cost of computing each of the $J{+}1$ splits is dominated by the cost of determining the Chernoff tau-leap step size, $\tau_{Ch}$ (see \eqref{eq:worktl}). As we observe in \cite{ourSL}, the work of computing a single $\tau_{Ch}$ is linear on $J$.
Then,  in order to avoid  $J^2$ complexity of the splitting rule, we implement a local search instead of computing $J$ $\tau_{Ch}$'s, to keep the complexity of $C_s$ linear on $J$.
The main idea is to keep track of the last split at each decision time, assuming that the propensities do not vary widely between.
If that is the case, we can just evaluate the previous split, $\mathcal{S}_\kappa$, and its neighbors, $\kappa{-}1$ and $\kappa{+}1$. Then, the cost of the splitting rule is on the order of three computations of a Chernoff step size. It turns out that this local search is very accurate for the examples we worked on. In order to avoid being trapped in local minima, a randomization rule may be applied.

\begin{rem}[Pareto Splitting rule]
\label{rem:pare}
Instead of computing a cost-based splitting at each decision time, the following rule can be applied:
$$\mathcal{R}_{\TL} \text{ is defined s.t. } \frac{\sum_{j\in \mathcal{R}_{\TL}} \tilde{a}_{\sigma(j)}}{\sum_{k=1}^J \tilde{a}_k} \geq \nu \COMMA$$
where $\nu$ is a problem-dependent threshold, which can be estimated using the cost-based splitting rule.
The idea is to use the tau-leap method for a $(100\times \nu) \%$  of the penalized activity (measured as before using the $\tilde{a}_j$'s), and an exact method for the other channels.

This rule is adaptive because it depends on the current state of the process, but it does not take into account the computational cost of the resulting partition of $\mathcal{R}$.
The advantage of this rule is that it is three times faster than the previous one. For the examples we worked on, the overall average gain in terms of computational work in a whole mixed path is about 45\% of the total work.

\end{rem}
\subsection{The one-step Mixing Rule}
In this section we present the main building block for simulating a mixed path. Let $x {=} \BX$ be the current state of the approximate process, $\bar X$. Therefore, the expected time step  of the MNRM is  given by $1/{a_0(x)}$. 
To move one step forward using the MNRM, we should compute at least $a_0(x)$ and sample a uniform random variable.
On the other hand, to move one step forward using the mixed Chernoff tau-leap method, we need first to compute the split, then compute the tau-leap increments for the reactions in the tau-leap set, $\mathcal{R}_\TL$, and finally compute the MNRM steps for the reactions in the set $\mathcal{R}_\MNRM$, as discussed in Section \ref{sec:cs}.

To avoid the overhead caused by unnecessary computation of the split, we first estimate the computational work of moving forward from the current time, $t$, to the next grid point, $\tilde{T}$, by using the MNRM only. If this work is less than the work of computing the split, we take an exact step. 
\begin{algorithm}[h!]
\caption{The one-step mixing rule. Inputs: the current state  of the approximate process, $\BX$, the current time, $t$, the values of the propensity functions evaluated at $\BX$, $\seqof{a_j(\BX)}{j=1}{J}$, the one-step exit probability bound $\delta$, the next grid point, $\tilde{T}$, and the previous optimal split, $\kappa$. Outputs: the tau-leap set, $\mathcal{R}_\TL$, the exact set, $\mathcal{R}_\MNRM$, and the new optimal split $\kappa$. 
}
\label{alg:sel}
\begin{algorithmic}[1]
	\REQUIRE $a_0 \leftarrow \sum_{j=1}^J a_j > 0$
	\IF{$K_1/a_0 < \tilde{T}-t$} 
		\STATE Compute $\theta_{j}$, {\small $j{=}1,..,J$} (see \eqref{eq:thetas})
		\STATE $\tilde{a}_{\sigma(j)} \leftarrow \text{Sort} \{(1{-}\theta_{j})a_j\}$ descending, {\small$j{=}1,..,J$}
		\STATE $\mathcal{S}_i \leftarrow $ Compute the splits taking into account the previous optimal split $\kappa$
		\STATE $(\mathcal{R}_\TL, \mathcal{R}_\MNRM, \kappa) \leftarrow $ Take the minimum work split
		\RETURN $(\mathcal{R}_\TL, \mathcal{R}_\MNRM,\kappa)$
	\ELSE
		\RETURN $(\emptyset, \mathcal{R},\kappa)$
	\ENDIF
\end{algorithmic}
\end{algorithm}
In order to compare the mentioned computational costs, we define $K_1$ as the ratio between the cost of computing the split, $C_s$, and the cost of computing one step using the MNRM.
\begin{rem}[Comparison with the one-step hybrid rule]
In \cite{ourSL} we developed a hybrid method, which, at each decision point, determines which method, exact or tau-leap, is cheaper to apply to the whole set of reactions. That is, in the hybrid method, we have either $\mathcal{R}_\TL =\emptyset$ and $\mathcal{R}_\MNRM =\mathcal{R}$ or $\mathcal{R}_\TL =\mathcal{R}$ and $\mathcal{R}_\MNRM =\emptyset$. Then, the mixed method can be seen as a generalization of the hybrid one. The key difference is in the cost of the decision rule, which, as we saw in Section \ref{sec:cs}, in the mixed method is on the order of three times the computation of the Chernoff step size. This difference can be significant in some problems. A Pareto splitting rule may be able to recover the cost of the hybrid one-step decision rule.
\end{rem}

\subsection{The Mixed-Path Algorithm}
In this section, we present a novel algorithm (Algorithm \ref{alg:mixed}) that combines the approximate Chernoff tau-leap method and the exact MNRM to generate a whole hybrid path. This algorithm combines the advantages of an exact method (expensive but exact) and the tau-leap method (may be cheaper but has a discretization error and a positive probability of exiting the lattice). 
This algorithm automatically and adaptively partitions the reactions into two  subsets, $\mathcal{R}_\TL$ and $\mathcal{R}_\MNRM$, using a computational work criterion.
Since a mixed path consists of a certain number of exact/approximate steps, it may also exit the lattice, except in those steps in which the tau-leap method is not applied; that is, when $\mathcal{R}_\TL$ is empty. 
The idea of this algorithm is to apply, at each decision point, the one-step mixing rule (Algorithm \ref{alg:sel}) to determine the sets $\mathcal{R}_\TL$ and $\mathcal{R}_\MNRM$, and then to apply the corresponding method.

\begin{algorithm}[h!]
\caption{The mixed-path algorithm. 
Inputs: the initial state, $X(0)$, the propensity functions, $\seqof{a_j}{j=1}{J}$, the stoichiometric vectors, $\nu{=}\seqof{\nu_j}{j=1}{J}$, the final time, $T$, and the one-step exit probability bound, $\delta$. 
Outputs: a sequence of states, $\seqof{\bar{X}(t_k)}{k=0}{K}$, 
and the number of times, $\NTL$, that the tau-leap method was successfully applied
(\ie, $\bar{X}(t_k)\in \latt$, we applied the tau-leap method and we obtained an $\bar{X}(t_{k+1})\in \latt$). 
Notes: given the current state, $next_{\MNRM}$ computes the next state using the MNRM method. Here, $t_i$ denotes the current time at the $i$-th step, and $\tau_{Ch}(\mathcal{R}_\TL)$ is the Chernoff step size associated with $\mathcal{R}_\TL$.
}
\label{alg:mixed}
\begin{algorithmic}[1]
\STATE $i \leftarrow 0$, $t_i \leftarrow t_0, \bar{X}(t_i) \leftarrow X(0), \bar Z \leftarrow X(0)$
\STATE $\mathcal{S}_j \leftarrow$ Compute splits, {\small $j{=}0,...,J$} 
\STATE $\kappa \leftarrow \arg \min_j \text{Work}(\mathcal{S}_j)$
\WHILE {$t_i<T$}
\STATE $\tilde T \leftarrow $ next grid point greater than $t_i$
\STATE $(\mathcal{R}_{\TL},\mathcal{R}_{\MNRM},\kappa) \leftarrow $ Algorithm \ref{alg:sel} with $(\bar Z,t_i,\seqof{a_j(\bar Z)}{j=1}{J},\delta, \tilde T,\kappa)$
\IF {$\mathcal{R}_{\TL}\neq \emptyset$}
\STATE $\Delta_{\TL} \leftarrow \mathcal{P}(a_j(\bar Z) \tau_{Ch}(\mathcal{R}_\TL)) \nu_{j}$,   for {\small$j{\in}\mathcal{R}_{\TL}$}
\STATE $H \leftarrow t_i+ \tau_{Ch}(\mathcal{R}_\TL)$
\ELSE
\STATE $H \leftarrow \min \{t_i{-}\log(r)/\sum_{j}a_j,T\}$,   $r{\sim} \text{Unif}(0,1)$
\ENDIF
\IF {$\mathcal{R}_{\MNRM}\neq \emptyset$}
\WHILE {$t_i<H$}
\STATE $(\bar{Z},t_i) \leftarrow next_{\MNRM}(\bar{Z},\mathcal{R}_e,t_i,H)$
\ENDWHILE
\ENDIF
\STATE $\bar{Z} \leftarrow \bar{Z}+\Delta_{\TL}$
\IF {$\bar Z \in \latt$}
	\STATE $\NTL \leftarrow \NTL+1$
	\STATE $t_{i+1} \leftarrow H$
\ELSE
	\RETURN $(\seqof{\bar{X}(t_k)}{k=0}{i},\NTL)$
\ENDIF

\STATE $i \leftarrow i+1$
\STATE $\bar X(t_i) \leftarrow \bar Z$
\ENDWHILE
\RETURN $(\seqof{\bar{X}(t_k)}{k=0}{i},\NTL)$
\end{algorithmic}
\end{algorithm}

\subsection{Coupled Mixed Paths}\label{cmp}
In this section, we explain how to couple two mixed paths. This is essential for the multilevel estimator. The four algorithms that are the building blocks of the coupling algorithm were already presented in \cite{ourML}. The novelty here comes from the fact that the coupled mixed algorithm may have to run the four algorithms concurrently in the sense of the time of the process, $t$. In this section, we denote with a bar $\bar \cdot$ and a double bar $\dbar \cdot$ coarse and fine grid-related quantities.

We now briefly describe the mixed Chernoff coupling algorithm, \ie,  Algorithm \ref{alg:coupled}. Let $\bar X$ and $\bar{\bar X}$ be two mixed paths, corresponding to two nested time discretizations, called coarse and fine, respectively.
Assume that the current time is $t$, and we know the states, $\bar X(t)$ and $\bar{\bar X}(t)$, the next grid points at each level, $\bar t$, $\dbar t$, and the corresponding one-step exit probabilities, $\bar \delta$ and $\dbar \delta$.
Based on this knowledge, we have to determine the four sets $(\bar{\mathcal{R}}_{\TL},\bar{\mathcal{R}}_{\MNRM},\dbar{\mathcal{R}}_{\TL},\dbar{\mathcal{R}}_{\MNRM})$, that correspond to four algorithms, B1, B2, B3 and B4, that we use as building blocks. Table \ref{tab:algs} summarizes them.
\begin{table}[h]
\centering
\begin{tabular}{l|c|c}
  & $\bar{\mathcal{R}}_{\TL}$ & $\bar{\mathcal{R}}_{\MNRM}$ \\ \noalign{\smallskip} \hline\noalign{\smallskip} 
$\dbar{\mathcal{R}}_{\TL}$ &B1 & B2 \\
$\dbar{\mathcal{R}}_{\MNRM}$ &B3 & B4
\end{tabular}
\medskip
\caption{Building blocks for simulating two coupled mixed Chernoff tau-leap paths.  
Algorithms B1 and B2 are presented as Algorithms 2 and 3 in \cite{Anderson2012}. Algorithms B3 and B4 can be directly obtained from Algorithm B2 (see \cite{ourML}).}\label{tab:algs}
\end{table}
In order to do that, the algorithm computes, independently, the sets $\mathcal{R}_\TL$ and $\mathcal{R}_\MNRM$ for each level, and the time until the next decision is taken, $H$, using Algorithm \ref{alg:timehor}. 
Next, it computes concurrently the increments due to each one of the sets (storing the results in $\Delta \bar X$ and $\Delta \dbar X$ for the coarse and fine grid, respectively). We note that the only case in which we use a Poisson random variates generator for the tau-leap method is in Algorithm B1 (Algorithm \ref{alg:bb1}). 
For Algorithms B2, B3 and B4, the Poisson random variables are simulated by adding independent exponential random variables with the same rate, $\lambda$, until exceeding a given time final time, $T$. The only difference in the latter blocks are the time points at which the propensities, $a_j$, are computed. For B2, the coarse propensities are frozen at time $t$, whereas for B3 the finer are frozen at $t$. In B4, the propensities are computed at each time step.
After arriving at time $H$, the four sets $(\bar{\mathcal{R}}_{\TL},\bar{\mathcal{R}}_{\MNRM},\dbar{\mathcal{R}}_{\TL},\dbar{\mathcal{R}}_{\MNRM})$ and the time until the next decision is taken, $H$, are determined again, and then all procedures are repeated until the simulation reaches the final time, $T$.

\newcommand{\SWITCH}[1]{\STATE \textbf{switch} #1 \begin{ALC@g}}
\newcommand{\ENDSWITCH}{\end{ALC@g}\STATE \textbf{end switch} }
\newcommand{\CASE}[1]{\STATE \textbf{case} #1\textbf{:} \begin{ALC@g}}
\newcommand{\ENDCASE}{\end{ALC@g}}

\begin{algorithm}[h!]
\caption{\small Coupled mixed path. Inputs: the initial state, $X(0)$, the final time $T$, the propensity functions, $\seqof{a_j}{j=1}{J}$, the stoichiometric vectors, $\seqof{\nu_j}{j=1}{J}$, and  
two time meshes, one coarser $\seqof{t_i}{i=0}{N}$, such that $t_N {=} T$ and a finer one, 
$\seqof{s_j}{j=0}{N'}$, such that $s_0 {=} t_0$, $s_M {=} t_N$, and $\seqof{t_i}{i=0}{N} {\subset} \seqof{s_j}{j=0}{N'}$.
Outputs: a sequence of states evaluated at the coarse grid, $\seqof{\bar{X}(t_k)}{k=0}{K} \subset \latt$, 
such that $t_K \leq T$, 
a sequence of states evaluated at the fine grid $\seqof{\dbar{X}(s_l)}{l=0}{K'} \subset \latt$,
such that $\bar{X}(t_K) \in \latt$ or $\dbar{X}(s_{K'}) \in \latt$.
If $t_K<T$, both paths exit the $\latt$ lattice before the final time, $T$. 
It also returns the number of times the tau-leap method was successfully applied at the fine level and at the coarse level and the number of exact steps at the fine level and at the coarse level.
For the sake of simplicity, we omit sentences involving the recording of current state variables, counting of the number of steps, checking if the path jumps out of the lattice, the updating of the current split, $\kappa$, and the return sentence.}
\label{alg:coupled}
\begin{algorithmic}[1]
\begin{small}
\STATE $t \leftarrow t_0$; $\bar{X} \leftarrow X(0)$; $\dbar{X} \leftarrow X(0)$ 
\STATE $\bar{t} \leftarrow$ next grid point in $\seqof{t_i}{i=0}{N}$ larger than $t$ 
\STATE $(\bar{H},\bar{\mathcal{R}}_{\TL},\bar{\mathcal{R}}_{\MNRM},\bar{a})  \leftarrow$ Algorithm \ref{alg:timehor} with ($\bar{X}$,$t$,$\bar{t}$,$T$,$\bar{\delta}$)
\STATE $\dbar{t} \leftarrow$ next grid point in $\seqof{s_i}{i=0}{N}$ larger than $t$ 
\STATE $(\dbar{H},\dbar{\mathcal{R}}_{\TL},\dbar{\mathcal{R}}_{\MNRM},\dbar{a})  \leftarrow$ Algorithm \ref{alg:timehor} with ($\dbar{X}$,$t$,$\dbar{t}$,$T$,$\dbar{\delta}$)

	\WHILE{$t<T$}
		\STATE $H \leftarrow \min\{\bar H, \dbar H\}$
		\STATE $(B_1,B_2,B_3,B_4) \leftarrow$ split building blocks from $(\bar{\mathcal{R}}_{\TL},\bar{\mathcal{R}}_{\MNRM},\dbar{\mathcal{R}}_{\TL},\dbar{\mathcal{R}}_{\MNRM})$
		\STATE Algorithm \ref{alg:bb1} (compute state changes due to block $B_1$)
		\STATE Initialize internal clocks $R, P$ if needed (see \cite{ourSL,ourML})
		\STATE $\Delta \bar{X} \leftarrow 0$; $\Delta \dbar{X} \leftarrow 0$
		\FOR {$\mathcal{B} = B_2,B_3,B_4$}
		\STATE $t_r \leftarrow t$
		\STATE $\bar{X}_r \leftarrow \bar{X}$; $\dbar{X}_r \leftarrow \dbar{X}$
		\WHILE{$t_r<H$}
		\STATE update $P_{j \in \mathcal{B}}$
		\SWITCH {$\mathcal{B}$}
			\CASE {$B_2$}
			 	\STATE $\bar{d} \leftarrow \bar{a}_{j\in \mathcal{B}}$
				\STATE $\dbar{d} \leftarrow a_{j\in \mathcal{B}}(\dbar{X})$
				\STATE $\tau_r \leftarrow $ Compute the Chernoff tau-leap step size using $(\bar{X}_r,\bar{a}_{j\in \mathcal{B}},H,\bar{\delta})$
			\ENDCASE
			\CASE {$B_3$}
				\STATE $\bar{d} \leftarrow a_{j\in \mathcal{B}}(\bar{X})$
				\STATE $\dbar{d} \leftarrow \dbar{a}_{j\in \mathcal{B}}$
				\STATE $\tau_r \leftarrow $ Compute the Chernoff tau-leap step size using $(\dbar{X}_r,\dbar{a}_{j\in \mathcal{B}},H,\dbar{\delta})$
			\ENDCASE
			\CASE {$B_4$}
				\STATE $\bar{d} \leftarrow a_{j\in \mathcal{B}}(\bar{X})$
				\STATE $\dbar{d} \leftarrow a_{j\in \mathcal{B}}(\dbar{X})$
				\STATE $\tau_r \leftarrow \infty$
			\ENDCASE
		\ENDSWITCH
			\STATE $A_1 \leftarrow \min(\bar{d},\dbar{d})$
			\STATE $A_2 \leftarrow \bar{d}-A_1$; $A_3 \leftarrow \dbar{d} -A_1$
			\STATE $H_r \leftarrow \min\{H,t_r{+}\tau_r\}$
				\STATE $(t_r, \bar{X}_r, \dbar{X}_r, R_{j \mathcal{B}}, P_{j\in \mathcal{B}}) \leftarrow$ Algorithm \ref{alg:auxilco} with $(t_r,H_r,\bar{X}_r,\dbar{X}_r,R_{j\in \mathcal{B}},P_{j\in \mathcal{B}},A)$
		\ENDWHILE
			\STATE $\Delta \bar{X} \leftarrow \Delta \bar{X} + (\bar{X}_r{-}\bar{X})$; $\Delta \dbar{X} \leftarrow \Delta \dbar{X} + (\dbar{X}_r{-}\dbar{X})$
		\ENDFOR
		\STATE $\bar{X} \leftarrow \bar{X} + \Delta \hat{X} + \Delta \bar{X}$;  $\dbar{X} \leftarrow \dbar{X} + \Delta \hat{\hat{X}} + \Delta \dbar{X}$
		\STATE $t \leftarrow H$
			\IF{$t<T$}
	\IF{$\bar H \leq \dbar H$}
		\STATE $\bar{t} \leftarrow$ next grid point in $\seqof{t_i}{i=0}{N}$ larger than $t$ 
\STATE $(\bar{H},\bar{\mathcal{R}}_{\TL},\bar{\mathcal{R}}_{\MNRM},\bar{a})  \leftarrow$ Algorithm \ref{alg:timehor} with ($\bar{X}$,$t$,$\bar{t}$,$T$,$\bar{\delta}$)

	\ENDIF
	\IF {$\bar H \geq \dbar H$}
		\STATE $\dbar{t} \leftarrow$ next grid point in $\seqof{s_j}{j=0}{N'}$ larger than $t$ 
\STATE $(\dbar{H},\dbar{\mathcal{R}}_{\TL},\dbar{\mathcal{R}}_{\MNRM},\dbar{a})  \leftarrow$ Algorithm \ref{alg:timehor} with ($\dbar{X}$,$t$,$\dbar{t}$,$T$,$\dbar{\delta}$)
	\ENDIF
\ENDIF
	\ENDWHILE
\end{small}
\end{algorithmic}
\end{algorithm}

\begin{algorithm}
\caption{\small Compute the next time horizon. Inputs: the current state, $\tilde{X}$, the current time, $t$, the next grid point, $\tilde{t}$, the final time, $T$, the one step exit probability bound, $\tilde{\delta}$, and the propensity functions, $a{=}\seqof{a_j}{j=1}{J}$.
Outputs: the next horizon $H$, the set of reaction channels to which the Tau-leap method should be applied, $\tilde{\mathcal{R}}_{\TL}$, the set of reaction channels to which MNRM should be applied, $\tilde{\mathcal{R}}_{\MNRM}$, and current propensity values $\tilde{a}$.}
\label{alg:timehor}
\begin{algorithmic}[1]
\begin{small}
\STATE $\tilde{a} \leftarrow a(\tilde{X})$
\STATE $(\tilde{\mathcal{R}}_{\TL},\tilde{\mathcal{R}}_{\MNRM}) \leftarrow$ Algorithm \ref{alg:sel} with $(\bar X,t,\seqof{a_j(\bar X)}{j=1}{J},\tilde \delta,\tilde t,\kappa)$
\IF{$\tilde{\mathcal{R}}_{\TL} \neq \emptyset$}
\STATE $\tilde{H} \leftarrow \min\{\tilde{t},t{+}\tau(\tilde{\mathcal{R}}_{\TL}),T\}$
\ELSE
\STATE $\tilde{H} \leftarrow \min\{t{+}\tau(\tilde{\mathcal{R}}_{\TL}),T\}$
\ENDIF
\RETURN $(\tilde H,\tilde{\mathcal{R}}_{\TL},\tilde{\mathcal{R}}_{\MNRM},\tilde a)$
\end{small}
\end{algorithmic}
\end{algorithm}

\begin{algorithm}[h!]
\caption{\small Compute building block 1. This algorithm is part of Algorithm \ref{alg:coupled}.
}
\label{alg:bb1}
\begin{algorithmic}[1]
\begin{small}
		\STATE $t_r \leftarrow t$
		\STATE $\Delta \hat{\hat{X}} \leftarrow 0$; $\Delta \hat{X} \leftarrow 0$
		\WHILE {$t_r < H$}
			\STATE $\bar{\tau}_r \leftarrow $ Compute the Chernoff tau-leap step size using $(\bar{X}{+}\Delta \hat{X},\bar{a}_{j\in B_1},H,\bar{\delta})$
			\STATE $\dbar{\tau}_r \leftarrow $ Compute the Chernoff tau-leap step size using $(\dbar{X}{+}\Delta \hat{\hat{X}},\dbar{a}_{j \in B_1},H,\dbar{\delta})$
			\STATE $H_r \leftarrow \min\{H,t_r{+}\bar{\tau}_r,t_r{+}\dbar{\tau}_r\}$
			\STATE $A_1 \leftarrow \min(\bar{a}_{j \in B_1},\dbar{a}_{j \in B_1})$
			\STATE $A_2 \leftarrow \bar{a}_{j \in B_1}-A_1$
			\STATE $A_3 \leftarrow \dbar{a}_{j \in B_1} -A_1$
			\STATE $\Lambda \leftarrow \mathcal{P}(A {\cdot} (H_r{-}t_r))$
			\STATE $\Delta \hat{\hat{X}} \leftarrow \Delta \hat{\hat{X}}+ (\Lambda_1 {+} \Lambda_2)\nu_{j \in B_1}$
			\STATE $\Delta \hat{X} \leftarrow \Delta \hat{X}+ (\Lambda_1 {+} \Lambda_3)\nu_{j \in B_1}$
			\STATE $t_r \leftarrow H_r$
		\ENDWHILE
\end{small}
\end{algorithmic}
\end{algorithm}

\begin{algorithm}[h!]
\caption{\small The auxiliary function used in algorithm \ref{alg:coupled}. 
Inputs: current time, $t$, current time horizon, $\dbar{T}$, current system state at coarser level and finer level, $\bar{X}$, $\dbar{X}$, respectively, the internal clocks $R$ and $P$, the values, $A$, and the current building block, $B$. Outputs: updated time, $t$, updated system states, $\bar{X}$, $\dbar{X}$, and updated internal clocks, $R_i$, $P_i$, $i{=}1,2,3$. 
}
\label{alg:auxilco}
\begin{algorithmic}[1]
\begin{small}
\STATE $\Delta t_i \leftarrow (P_i - R_i)/A_i$, for $i=1,2,3$
\STATE $\Delta \leftarrow \min_i\{ \Delta t_i \}$
\STATE ${\mu} \leftarrow \text{argmin}_i \{\Delta t_i\}$
					\IF{$t+\Delta > \dbar{T}$}
						\STATE $R \leftarrow R + A {\cdot} (\dbar{T}{-}t)$
						\STATE $t \leftarrow \dbar{T}$
					\ELSE
						\STATE update $\bar{X}$ and $\dbar{X}$ using $\nu_{j \in B}$
						\STATE $R \leftarrow R + A \Delta$
						\STATE $r \leftarrow$ uniform$(0,1)$
						\STATE $P_\mu \leftarrow P_\mu + \log (1/r)$
						\STATE $t \leftarrow t + \Delta$
					\ENDIF
					\RETURN $(t, \bar X, \dbar X, R, P)$
\end{small}
\end{algorithmic}
\end{algorithm}


\section{The Multilevel Estimator and Total Error Decomposition}
\label{sec:mlmcerr}
In this section, we first show the multilevel Monte Carlo estimator. We then  analyze and control the computational global error, which is decomposed into three error components: the discretization error, the global exit error, and the Monte Carlo statistical error. Upper bounds for each one of the three components are given. 
Finally, we briefly describe the automatic estimation procedure that allows us to estimate our quantity of interest within a given prescribed relative tolerance, up to a given confidence level.


\subsection{The MLMC Estimator}\label{sec:mlmc}
In this section, we discuss and implement a multilevel Monte Carlo estimator for the mixed Chernoff tau-leap case. 
Consider a hierarchy of nested meshes of the time interval $[0,T]$, indexed by $\ell = 0,1,\ldots,L$.
Let $\Delta t_0$ be the size of the coarsest time mesh that corresponds to the level $\ell{=}0$.
The size of the time mesh at  level $\ell \geq 1$ is given by $\Delta t_{\ell} {=} R^{-\ell} \Delta t_0$, where $R{>}1$ is a given integer constant. 
Let $\{\bar X_{\ell}(t)\}_{t\in[0,T]}$ be a mixed Chernoff tau-leap process with a time mesh of size $\Delta t_{\ell}$ and a one-step exit probability bound $\delta$, and let $g_{\ell}{:=} g({\bar X}_{\ell}(T))$ be our quantity of interest computed with a mesh of size $\Delta t_{\ell}$. We can simulate paths of  $\{\bar X_{\ell}(t)\}_{t\in[0,T]}$ by using  Algorithm \ref{alg:mixed}.
We are interested in estimating $\expt{g_{L}}$, and we can simulate correlated pairs, $(g_{\ell},g_{\ell{-}1})$ for $\ell=1,\ldots,L$, by using Algorithm \ref{alg:coupled}.
Let $A_\ell$ be the event in which the $\ell$-th grid level path, $\bar X_{{\ell}}$, arrives at the final time, $T$, without exiting the state space of $X$. 

Consider the following telescopic decomposition:
\begin{align*}
\expt{g_L \indicator{A_L}} = \expt{g_0 \indicator{A_0}} {+} \sum_{\ell=1}^L \expt{g_\ell\indicator{A_\ell}-g_{\ell-1}\indicator{A_{\ell-1}}} \COMMA
\end{align*} 
where $\indicator{A}$ is the indicator function of the set $A$. This motivates the definition of our MLMC estimator of $\expt{g(X(T))}$:
\begin{align}\label{MLMCest}
\mathcal{M}_L := \frac{1}{M_0} \sum_{m=1}^{M_0}g_0\indicator{A_0}(\omega_{m,0}) + \sum_{\ell=1}^L \frac{1}{M_{\ell}} \sum_{m=1}^{M_{\ell}} [g_{\ell} \indicator{A_\ell} - g_{\ell-1} \indicator{A_{\ell-1}}](\omega_{m,\ell}) \PERIOD
\end{align}
\subsubsection*{Computational Complexity} 
\label{sec:complex}
A key property of our multilevel estimator is that the computational work is a function of the given relative tolerance, $TOL$, is of the order of $TOL^{-2}$.
 The optimal work is is given by 
\begin{equation*}
w^*_{L} (TOL) = \LP\frac{C_A}{\theta}\sum_{\ell=0}^L \sqrt{\mV{\ell} \psi_{\ell}}\RP^{2} TOL^{-2}\PERIOD
\end{equation*}
From the fact that the sum $\sum_{\ell=0}^{\infty} \sqrt{\mV{\ell} \psi_{\ell}}$ converges, because $\psi_\ell = \Ordo{\psi_{\MNRM}}$, 
we conclude that 
$\sup_{L}\{\sum_{\ell=0}^L \sqrt{\mV{\ell} \psi_{\ell}}\}$ is bounded and, therefore, the expected computational complexity of the multilevel mixed Chernoff tau-leap method is $w^*_L(TOL){=}\Ordo{TOL^{-2}}$. 

\subsection{Global Error Decomposition}
In this section, we define the computational global error, $\EE_L$, and show how it can be naturally decomposed into three components: the discretization error, $\EE_{I,L}$, and the exit error, $\EE_{E,L}$, both coming from the tau-leap part of the mixed method, and the Monte Carlo statistical error, $\EE_{S,L}$. 
We also give upper bounds for each one of the three components. 

The computational global error, $\EE_L$, is defined as
\begin{equation*}
\EE_L := \expt{g(X(T))} - \mathcal{M}_L \COMMA
\end{equation*}
and can be decomposed as
\begin{align*}
\expt{g(X(T))} - \mathcal{M}_L &= \expt{g(X(T))(\indicator{A_L}+\indicator{A^c_L})} \pm \expt{g_L \indicator{A_L}} - \mathcal{M}_L\nonumber \\
&= \underbrace{\expt{g(X(T)) \indicator{A^c_L}}}_{=:\EE_{E,L}} 
 + \underbrace{\expt{ \LP g(X(T)) {-} g_L \RP \indicator{A_L}}}_{=:\EE_{I,L}} + \underbrace{\expt{g_L \indicator{A_L}}{-}\mathcal{M}_L}_{=:\EE_{S,L}}\PERIOD 
\end{align*}

We showed in \cite{ourSL} that by choosing adequately the one-step exit probability bound, $\delta$, the exit error, $\EE_{E,L}$, satisfies 
$|\EE_{E,L}| \leq |\expt{g(X(T))}|\, \prob{A^c_L}\leq TOL^2$. 

An efficient procedure for accurately estimating $\EE_{I,L}$ in the context of the tau-leap method is described in \cite{ourML}. 
For each mixed path,  
$\seqof{\bar X_{\ell}(t_{n,\ell},\bar \omega)}{n=0}{N(\bar{\omega})}$, we define  the sequence of dual weights, $\seqof{\varphi_{n,\ell}(\bar \omega)}{n=1}{N(\bar{\omega})}$, backwards as follows:
\begin{align}\label{eq:phisback}
\varphi_{N(\bar{\omega}),\ell} &:= \nabla g(\bar X_{\ell}(t_{N(\bar{\omega}),\ell},\bar \omega))\\
\varphi_{n,\ell} &:= \LP Id +  {\Delta t_{n,\ell}} \JAC_a^T(\bar X_{\ell}(t_{n,\ell},\bar \omega)) \,\nu^T\RP \varphi_{n+1,\ell}, \quad n=N(\bar{\omega}){-}1,\ldots,1\COMMA\nonumber
\end{align}
where $\Delta t_{n,\ell}{:=}t_{n+1,\ell}{-}t_{n,\ell}$, $\nabla$ is the gradient operator and $\JAC_a(\bar X_{\ell}(t_{n,\ell},\bar \omega)){\equiv}[\partial_i a_j(\bar X_{\ell}(t_{n,\ell},\bar \omega))]_{j,i}$ is the Jacobian matrix of the propensity function, $a_j$, for $j{=}1\ldots J$ and $i{=}1\ldots d$. 
We then approximate $\EE_{I,L}$ by $\avg{\EE_{I,L}(\bar\omega)}{\cdot}$, where 
\begin{equation*}
\EE_{I,L}(\bar \omega) := \sum_{n=1}^{{N(\bar{\omega})}}  \LP\frac {\Delta t_{n,L}}{2} \varphi_{n,L}   \sj \indicator{j\in \mathcal{R}_{\TL}(n)} \nu_j^T \LP a_j(\bar X_{L}(t_{n+1,\ell})){-}a_j(\bar X_{L}(t_{n,\ell}))\RP \RP(\bar \omega)\COMMA
\end{equation*}
$\avg{X}{M}{:=}\frac 1 M \sum_{m=1}^M X(\omega_m)$ and $\svar{X}{M}{:=}\avg{X^2}{M}- \avg{X}{M}^2$ denote
the sample mean and the sample variance of the random  variable, $X$, respectively.
Here $\indicator{j\in \mathcal{R}_{\TL}(n)}{=}1$ if and only if, at time $t_{n,\ell}$, the tau-leap method was used for reaction channel $j$, and we denote by $Id$ the $d\times d$ identity matrix. 

The variance of the statistical error, $\EE_{S,L}$, is given by $\sum_{\ell=0}^L \frac{\mV{\ell}}{M_{\ell}}$,
where $\mV{0}:=\var{g_{0}\indicator{A_0}}$ and 
$\mV{\ell}:=\var{ g_{\ell} \indicator{A_{\ell}} -g_{\ell-1} \indicator{A_{\ell{-}1}}},\,\,\ell \geq 1$. 
In \cite{ourML}, we presented an efficient and accurate method for estimating $\mV{\ell}$, $\ell \geq 1$ using the formula
\begin{align*}
\hat{\mathcal{V}}_\ell&:=\svar{\sum_n \expt{ \varphi_{n+1} \cdot e_{n+1}\SEP\mathcal{F}}(\bar \omega)}{M_{\ell}}+ \avg{\sum_n \var{ \varphi_{n+1} \cdot e_{n+1}\SEP\mathcal{F}}(\bar \omega)}{M_{\ell}}\COMMA
\end{align*}
where $\mathcal{F}$ is a suitable chosen sigma algebra such that $\seqof{\varphi_n(\bar \omega)}{n=1}{N(\bar \omega)}$ is measurable, with $N(\bar \omega)$ being the total number of steps given by Algorithm \ref{alg:coupled}. In this way, the only randomness in $\expt{ \varphi_{n+1} \cdot e_{n+1}\SEP\mathcal{F}}$ and $\var{ \varphi_{n+1} \cdot e_{n+1}\SEP\mathcal{F}}$ comes from the local errors, $\seqof{e_n}{n=1}{N(\bar \omega)}$, defined as $e_{n} := X_{\ell,n} - X_{\ell-1,n}$.
In the aforementioned work, we derived exact and approximate formulas for computing $\expt{ \varphi_{n+1} \cdot e_{n+1}\SEP\mathcal{F}}$ and $\var{ \varphi_{n+1} \cdot e_{n+1}\SEP\mathcal{F}}$.
\begin{rem}[Backward Euler] In \eqref{eq:phisback}, we have that $\varphi_{n,\ell}$ can be computed by a backward Euler formula when too fine time meshes are required for stability, \ie,
$\varphi_{n,\ell} := \LP Id -  {\Delta t_{n,\ell}} \JAC_a^T(\bar X_{\ell}(t_{n,\ell},\bar \omega)) \,\nu^T\RP^{-1} \varphi_{n+1,\ell}$.   
\end{rem}

\subsection{Estimation Procedure}
\label{sec:estim}
In this section, we briefly describe the automatic procedure that estimates $\expt{g(X(T))}$ within a given prescribed relative tolerance, $TOL{>}0$, up to a given confidence level. Up to minor changes, it is the same as the one presented in \cite{ourML}. It is important to remark that the minimal user intervention is required to obtain the parameters needed to simulate the mixed paths, and subsequently, to compute the estimations using \eqref{MLMCest}. 
Once the reaction network is given (stoichiometric matrix $\nu$ and $J$ propensity functions $a_j$), the user only needs to set the required maximum allowed relative global error or tolerance, $TOL$, and the confidence level, $\alpha$. This process has three phases:
\begin{description}
\item[Phase I] Calibration of virtual machine-dependent quantities. In this phase, we estimate the quantities 
$C_\MNRM$, $C_\TL$, $C_s$ and the function $C_P$
that allow us to model the expected computational work, measured in runtime. 
\item[Phase II] Solution of the work optimization problem:  we obtain 
the total number of levels, $L$, and the sequences, $\seqof{\delta_\ell}{\ell=0}{L}$ and $\seqof{M_\ell}{\ell=0}{L}$, \ie, the one-step exit probability bounds and the required number of simulations at each level. 
In this phase, given a relative tolerance, $TOL{>}0$, we solve the work optimization problem
\begin{align}
\left\{ \!
\begin{array}{l}\label{eq:minworkprob}
\min_{\{\Delta t_0, L,\seqof{M_{\ell},\delta_{\ell}}{\ell=0}{L}\}} \sum_{\ell =0}^{L}\psi_\ell M_{\ell}\\
\mbox{s.t.}\\
\EE_{E,L} +\EE_{I,L}+\EE_{S,L} \leq TOL
\end{array}
\right. \PERIOD
\end{align}
An algorithm to efficiently compute the solution of this optimization problem is given in \cite{ourML}.
Our objective function is the expected total work of the MLMC estimator, $\mathcal{M}_L$, \ie,
$\sum_{\ell =0}^{L}\psi_\ell M_{\ell}\COMMA$
where $L$ is the deepest level, $\psi_0$ is the expected work of a single-level path at level 0, and $\psi_\ell$, for $\ell\geq 1$, is the expected computational work of two coupled paths at levels $\ell{-1}$ and $\ell$. Finally, $M_0$ is the number of single-level paths at level 0, and $M_\ell$, for $\ell\geq 1$, is the number of coupled paths at levels $\ell{-1}$ and $\ell$. 
We now describe the quantities $\seqof{\psi_\ell}{\ell=0}{L}$. 
First, $\psi_0$ is the expected work of a single hybrid path (simulated by Algorithm \ref{alg:mixed}), 
\begin{align}\label{eq:cost}
\psi_0 &:= C_\MNRM \expt{\NMNRM(\Delta t_0,\delta_0)} + C_\TL \expt{\NTL (\Delta t_0,\delta_0)}\\ &+  \int_{[0,T]}\expt{ \sum_{j \in \mathcal{R}_{\TL}(s)} C_P(a_j(\bar X_0(s))\tau_{Ch}(\bar X_0(s),\delta_0)) ds} \nonumber
\COMMA
\end{align} 
where $\Delta t_0$ is the size of the time mesh at level 0 and $\delta_0$ is the exit probability bound at level 0, and $\mathcal{R}_\TL=\mathcal{R}_\TL(t)$ is the tau-leap set, which depends on time (and also the current state of the process).
The set $\mathcal{R}_\TL$ is determined  at each decision step by Algorithm \ref{alg:sel}.
Therefore, the expected work at level 0 is $\psi_0 M_0$, 
where $M_0$ is the total number of single hybrid paths. 

For $\ell \geq 1$, we use Algorithm \ref{alg:coupled} to generate $M_{\ell}$-coupled paths that couple levels $\ell{-}1$ and $\ell$.
The expected work of a pair of coupled hybrid paths at levels $\ell$ and $\ell-1$ is
\begin{align}\label{eq:costc}
\psi_\ell &:= C_\MNRM \expt{\NMNRMC(\ell)} +C_\TL \expt{\NTLC(\ell)} \\ \nonumber
&+  \int_{[0,T]} \expt{\sum_{j\in \mathcal{R}_{\TL,\ell}(s)} C_P(a_j(\bar X_\ell(s))\tau_{Ch}(\bar X_\ell(s),\delta_\ell)) ds}\\ \nonumber
 &+ \int_{[0,T]} \expt{\sum_{j\in \mathcal{R}_{\TL,\ell-1}(s)}  C_P(a_j(\bar X_{\ell-1}(s))\tau_{Ch}(\bar X_{\ell-1}(s),\delta_{\ell-1})) ds}\COMMA
\end{align} 
where 
\begin{align*}
&\NMNRMC(\ell):=\NMNRM(\Delta t_\ell,\delta_\ell) + \NMNRM(\Delta t_{\ell-1},\delta_{\ell-1})\\
&\NTLC(\ell) :=\NTL(\Delta t_\ell,\delta_\ell)+\NTL (\Delta t_{\ell-1},\delta_{\ell-1})\PERIOD
\end{align*} 
\item[Phase III] Estimation of $\expt{g(X(T))}$. 
\end{description} 

\section{A Control Variate Based on a Deterministic Time Change}
\label{sec:cvar}
In this section, we motivate a novel control variate for the random variable $X(T,\omega)$ defined by the random time change representation,
\begin{equation*}
X(T,\omega) = x_0 + \sum_j \nu_j Y_j\LP \int_0^T a_j(X(s))\,ds,\omega\RP\PERIOD
\end{equation*}

First, we replace the independent Poisson processes, $\LP Y_j(s,\omega) \RP_{s\geq 0}$, by the identity function. This defines the deterministic mean field, 
\begin{equation*}
Z(T) = x_0 + \sum_j \nu_j  \int_0^T a_j(Z(s))\,ds\PERIOD
\end{equation*}

Next,  we consider the random variable 
\begin{equation*}
\tilde X(T,\omega) = x_0 + \sum_j \nu_j Y_j\LP \int_0^T a_j(Z(s))\,ds, \omega\RP\COMMA
\end{equation*}
which uses the same realizations of $\LP Y_j(s,\omega) \RP_{s\geq 0}$ that define $X(T,\omega)$. 
In this way, we expect some correlation between  $X(T)$ and $\tilde X(T)$. 
Since $\expt{\tilde X(T)} = Z(T)$ is a computable quantity, we have that $\tilde X(T)$ is a potential control variate for $X(T)$ obtained at almost negligible extra computational cost.

We have that $\tilde X(T,\omega)$ can be considered as a deterministic time change approximation of $X(T,\omega)$.

To implement this idea, we first consider the sequence $Z_{k}$, defined as a forward Euler discretization of the mean field over a suitable mesh,
 $\{t_0{=}0,  t_1,\ldots, t_K{=}T\}$, $\Delta t_k :=  t_{k+1}{-} t_k$, $k{=}0,1,\ldots,K{-}1$; that is,
$$\left\{ 
\begin{array}{ll}
Z_{k+1} = Z_{k} + \sum_j \nu_j a_j(Z_{k}) \Delta t_{k}, & k{=}0,\ldots,K{-}1\\
Z_0 = x_0  &
\end{array}
\right. .
$$
The sequence  $Z_{k}$ allow us to define another sequence, $\hat \Lambda_{j,k}$, by 
$$\left\{ 
\begin{array}{ll}
\hat \Lambda_{j,k+1} = \hat \Lambda_{j,k} + a_j(Z_k)\Delta t_{k}, & k{=}1,\ldots,K{-}1\\
\hat \Lambda_{j,0} = 0  &
\end{array}
\right. ,
$$
where $\hat \Lambda_{j,K}$ approximates $\int_0^T a_j(Z(s))\,ds$.

Then, for each realization of $\bar{X}(T,\omega)$, which is an approximation of ${X}(T,\omega)$, we compute the control variate:
\begin{equation}\label{eq:XK}
\hat X_K = x_0 + \sum_j \nu_j Y_j\LP \hat \Lambda_{j,K}\RP\COMMA
\end{equation}
which is the corresponding approximation of  $\tilde{X}(T,\omega)$ and  has the computable expectation 
$$\mu_K := \expt{\hat X_K} = x_0 + \sum_j \nu_j  \hat \Lambda_{j,K}. $$

Now, we consider the random sequence, $\{\bar X_n(\omega)\}_{n=0}^{N(\omega)}$, generated in this case by the mixed algorithm. Here, $\bar X(\omega)_{N(\omega)}$ is an approximation of $X(T,\omega)$.
The sequence of mixed random times, $\{ \bar \Lambda_{j,n}(\omega) \}$, is defined by
$$\left\{ 
\begin{array}{ll}
\bar \Lambda_{j,n+1} = \bar \Lambda_{j,n} + a_j(\bar X_n(\omega))\Delta s_{n}, & n{=}0,\ldots,N(\omega){-}1\\
\bar \Lambda_{j,0} = 0 &
\end{array}
\right. \COMMA
$$
over the mesh $\{s_0{=}0,  s_1,\ldots, s_{N(\omega)}{=}T\}$, $\Delta s_n :=  s_{n+1}{-} s_n$, $n{=}0,1,\ldots,N(\omega){-}1$.


At this point, it is crucial to observe that we can keep track of the values $Y_j(\bar \Lambda_{j,n},\omega)$, since at each step of the approximation algorithm, we are sampling the increments of the processes, $Y_j$. From now on, we omit $\omega$ in our notation.

The values $Y_j\LP \hat \Lambda_{j,K}\RP$, required in \eqref{eq:XK}, can by obtained by sampling the process $Y_j$ as follows.
For each realization of $\bar{X}$, we have two scenarios:
\begin{enumerate}
\item for some $n$, $\bar \Lambda_{j,n} <\hat \Lambda_{j,K}< \bar \Lambda_{j,n+1}$.  Since $Y_j\LP \bar \Lambda_{j,n} \RP$ and $Y_j\LP \bar \Lambda_{j,n+1} \RP$ are known, we sample a Poissonian bridge (binomial), \ie,
$$Y_j(\hat \Lambda_{j,K})\sim Y_j\LP \bar \Lambda_{j,n} \RP + \text{binomial}\LP Y_j\LP \bar \Lambda_{j,n+1} \RP-Y_j\LP \bar \Lambda_{j,n} \RP, \frac{\hat \Lambda_{j,K} - \bar \Lambda_{j,n}}{\bar \Lambda_{j,n+1} - \bar \Lambda_{j,n}}\RP. $$ 
\item $\bar \Lambda_{j,K} >\hat \Lambda_{j,N}$. Since we know the value $Y_j\LP \bar \Lambda_{j,N} \RP$,
we just have to sample a Poisson random variate as follows:
$$Y_j(\hat \Lambda_{j,K})\sim Y_j\LP \bar \Lambda_{j,N} \RP + \text{Poisson}(a_j(\bar X_N)(\hat \Lambda_{j,K}- \bar \Lambda_{j,N})).$$
\end{enumerate}


Finally, using the aforementioned control variate, we can estimate $\expt{g(\bar{X}(T))}$ with 
$$\frac{1}{M} \sum_{m=1}^{M} g(\bar{X}_N(\omega_m)) - \beta \frac{1}{M}\sum_{m=1}^{M} (g(\hat{X}_K(\omega_m))-g(\mu_K)) \COMMA$$
for any linear functional, $g$. 
For polynomial observables, $g$, this estimator can be easily extended by Taylor expansions in such way that we can estimate $\expt{g(\tilde X(T))}$ by powers, $g\LP \LP\expt{\tilde X(T)}\RP^k\RP$.

\begin{rem}[Reducing the variance at the coarsest level]
The main application of the deterministic time change control variate, $\tilde X(T)$, in this work is at the coarsest level of our multilevel hierarchy. 
Consider the trivial decomposition 
$$g(\bar X_0(T)) = g(\tilde X(T)) + \LP g(\bar X_0(T))-  g(\tilde X(T)) \RP.$$ Therefore, 
$$\expt{g(\bar X_0(T))} = \expt{g(\tilde X(T))} + \expt{g(\bar X_0(T))-  g(\tilde X(T))}.$$ 
Since we can compute exactly $ \expt{g(\tilde X(T))}$, we just have to estimate $\expt{g(\bar X_0(T))- g(\tilde X(T))}$ instead of $\expt{g(\bar X_0(T))}$ in our multilevel scheme. The computational gain lies in the fact that $\var{g(\bar X_0(T))- g(\tilde X(T))}$ could be substantially lower than $\var{g(\bar X_0(T))}$.
\end{rem}
\begin{rem}[Computational Cost]
An advantage of this control variate is that the computational cost is almost negligible because we only need to store two scalars, $\bar \Lambda_{j,n}$ and $\bar \Lambda_{j,n+1}$, for each reaction, $j$. These values are determined at each step by $a_j(\bar X_n)$, which is a quantity that is already computed at each time step of the mixed algorithm. Also, for each realization of the control variate, at most one Poisson random variate is needed for each reaction channel.
\end{rem}
\begin{rem}[Empirical Time Change]
We can also compute the final times, $\hat \Lambda_{j,K}$, using a sample average of mixed paths instead of the mean field. We found no significant improvements when using that approach, which requires a lot more computational work. We conjecture that, for settings in which the mean field is not representative, this approach is the only reasonable option.
\end{rem}

\section{Numerical Examples}
\label{sec:examples}

\newcommand{\MLVirus}{MLM_Virus_2014_4_28_10_14_33}
\newcommand{\Stiff}{MLM_Simple_stiff_2014_5_25_11_7_22}

In this section, we present two examples to illustrate the performance of our proposed method, and we compare the results with the hybrid MLMC approach given in \cite{ourML}. For benchmarking purposes, we use Gillespie's Stochastic Simulation Algorithm (SSA) instead of the Modificed Next Reaction Method (MNRM) because the former is widely used in the literature.


\subsection*{Intracellular Virus Kinetics}

This model, first developed in \cite{srivastava2002stochastic}, has four species and six reactions,
\begin{itemize}
\item $E \xrightarrow{1} E {+} G$, the viral template (E) forms a viral genome (G),
\item $G \xrightarrow{0.025} E$, the genome generates a new template,
\item $E \xrightarrow{1000} E {+} S$, a viral structural protein (S) is generated,
\item $G{+}S \xrightarrow{7.5\times10^{-6}} V$, the virus (V) is produced, 
\item $E\xrightarrow{0.25}\emptyset$, $S\xrightarrow{2}\emptyset$ degradation reactions.
\end{itemize}
Its stoichiometric matrix and its propensity functions, $a_j:\zset_+ \rightarrow \rset$, are given by 
\begin{align*}
\nu =\left( 
 \begin{array}{ccccc}     1    & 0 &    0 & 0\\
     -1    & 0 &    1 & 0 \\
     0   & 1 &    0 & 0 \\
     -1 & -1 &0 &1 \\
    0    & 0 &   -1&  0 \\
     0   & -1 &  0 & 0 
 \end{array} 
 \right) ^{tr}\mbox{    and     }  a(X) =\left( \begin{array}{c} E \\ 0.025 \, G \\ 1000 \, E \\ 7.5{\times} 10^{-6} G \, S \\ 0.25 \, E \\ 2 \, S \end{array} \right)\COMMA
\end{align*}
respectively. 

In this model, $X(t)=(G(t),S(t),E(t),V(t))$, and $g(X(t))=V(t)$, the number of viruses produced. The initial condition is $X_0{=}(0, 0, 10, 0)$ and the final time is $T{=}20$.
This example is interesting because i) it shows a clear separation of time scales, ii) our previous hybrid Chernoff method has no compuational work gain with respect to an exact method, and iii) in \cite{Anderson2012} the authors take an alternative approach, not using the multilevel aspect of their paper.

\begin{figure}[h!]
\centering
\begin{minipage}{0.49\textwidth}
\includegraphics[scale=0.33]{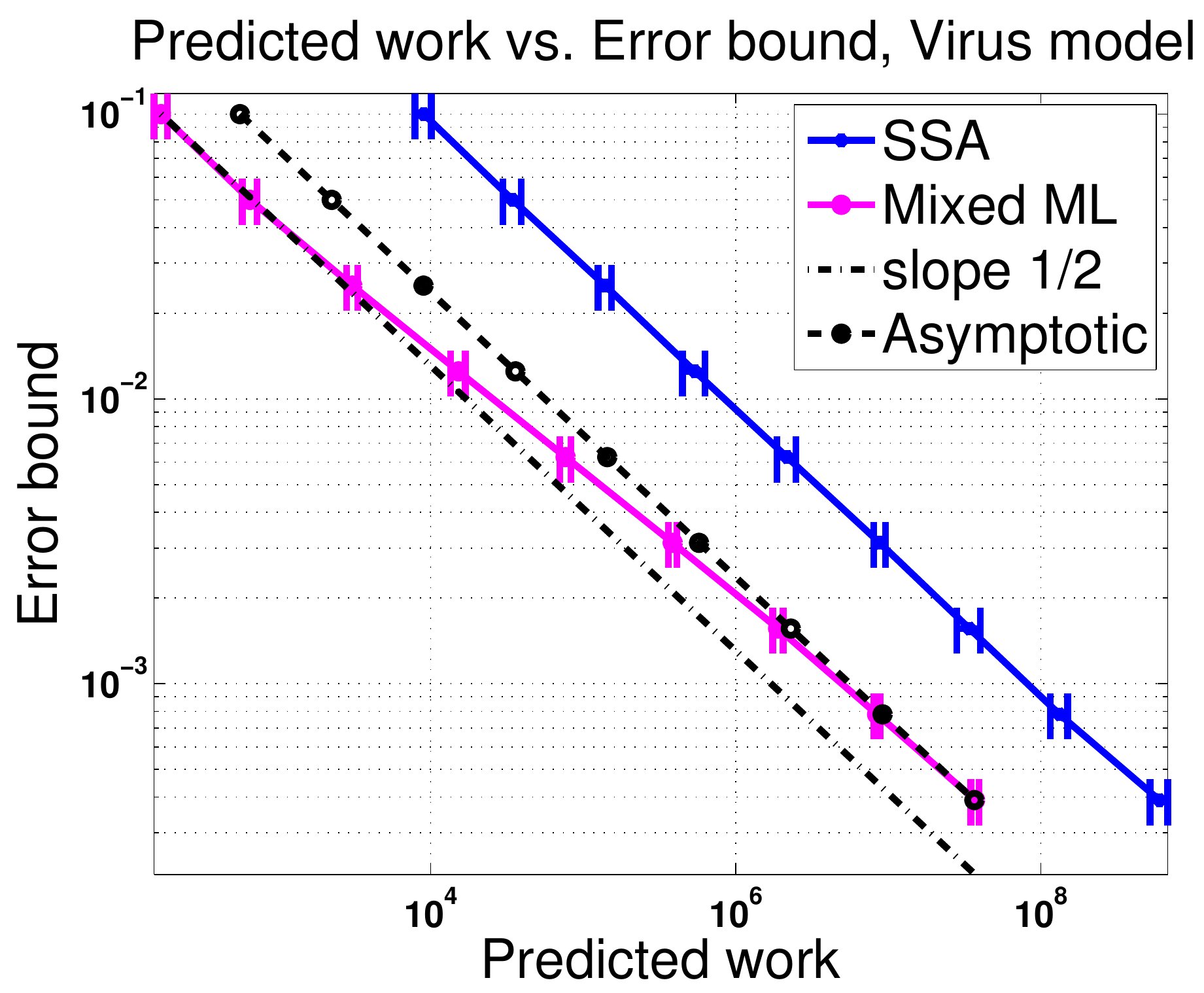}
\end{minipage}
\hfill
\begin{minipage}{0.49\textwidth}
\begin{tabular}{lcc}
\label{tab:vir}
$TOL$ & $L^*$ & $\frac{\hat{W}_{ML}}{\hat{W}_{\ssa}}$ \\ \noalign{\smallskip} \hline\noalign{\smallskip} 
1.00e-01 & 1.0  & 0.02 $\pm$0.001   \\ 
5.00e-02 & 1.0  & 0.02 $\pm$0.001   \\ 
2.50e-02 & 1.2 $\pm$0.261 & 0.02 $\pm$0.001 \\ 
1.25e-02 & 2.2 $\pm$0.261 & 0.03 $\pm$0.002   \\ 
6.25e-03 & 3.4 $\pm$0.320 & 0.04 $\pm$0.004  \\ 
3.13e-03 & 4.6 $\pm$0.320 & 0.04 $\pm$0.002   \\ 
1.56e-03 & 5.8 $\pm$0.261 & 0.06 $\pm$0.008  \\ 
7.81e-04 & 7.4 $\pm$0.433 & 0.07 $\pm$0.006  \\ 
3.91e-04 & 8.6 $\pm$0.320 & 0.06 $\pm$0.007   \\ 
\noalign{\smallskip}\hline 
\end{tabular}
\end{minipage}
\caption{Left: Predicted work (runtime) versus the estimated error bound, with $95\%$ confidence intervals. The multilevel mixed method is preferred over the SSA and the multilevel hybrid  method for all the tolerances. 
Right: Details for the ensemble run of the phase II algorithm.
Here, $\hat{W}_{ML}= \sum_{\ell =0}^{L^*}\hat \psi_\ell M_{\ell}$ and  $\hat{W}_{\ssa} = M_{\ssa}\, \,C_\ssa \,\avg{\NSSAP}{\cdot}$. 
As an example, the fourth row of the table tells us that, for a tolerance $TOL{=}1.25\cdot 10^{-2}$, 2.2 levels are needed on average. 
The work of the multilevel hybrid method is, on average, $3\%$ of the work of the SSA and the multilevel hybrid method. Confidence intervals at 95\% are also provided.
}
\label{fig:dec2-worktime}
\end{figure}

We now analyze an ensemble of 10 independent runs of the phase II algorithm (see Section \ref{sec:estim}), using different relative tolerances. 
In Figure \ref{fig:dec2-worktime}, we show the total predicted work (runtime) for the multilevel mixed method and for the SSA method, versus the estimated error bound. We also show the estimated asymptotic work of the multilevel mixed method. We remark that the computational work of the multilevel \emph{hybrid} method is the same as the work of the SSA. 

%

\begin{figure}[h!]
\centering
\begin{minipage}{0.49\textwidth}
\includegraphics[scale=0.31]{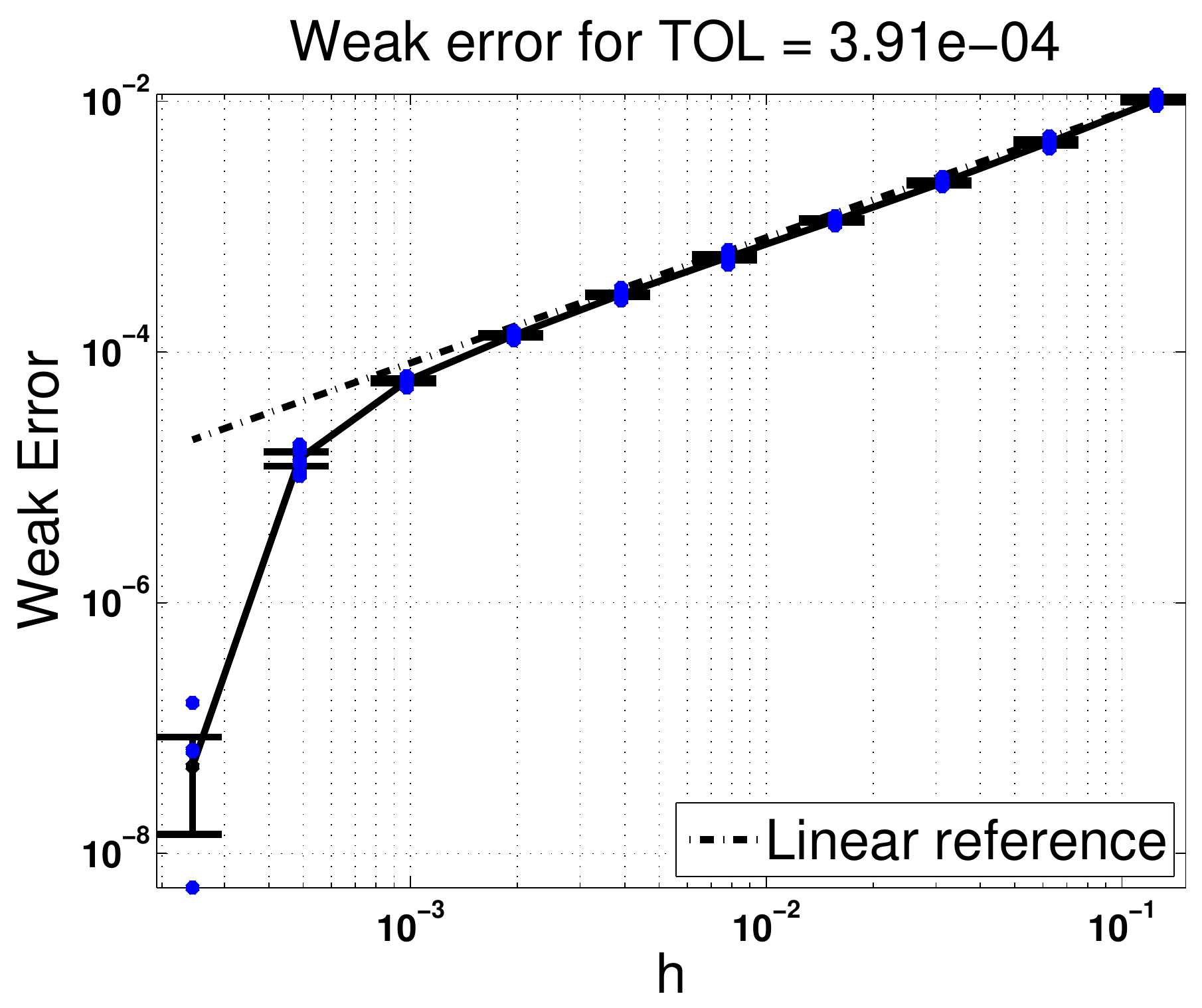}
\end{minipage}
\begin{minipage}{0.49\textwidth}
\includegraphics[scale=0.31]{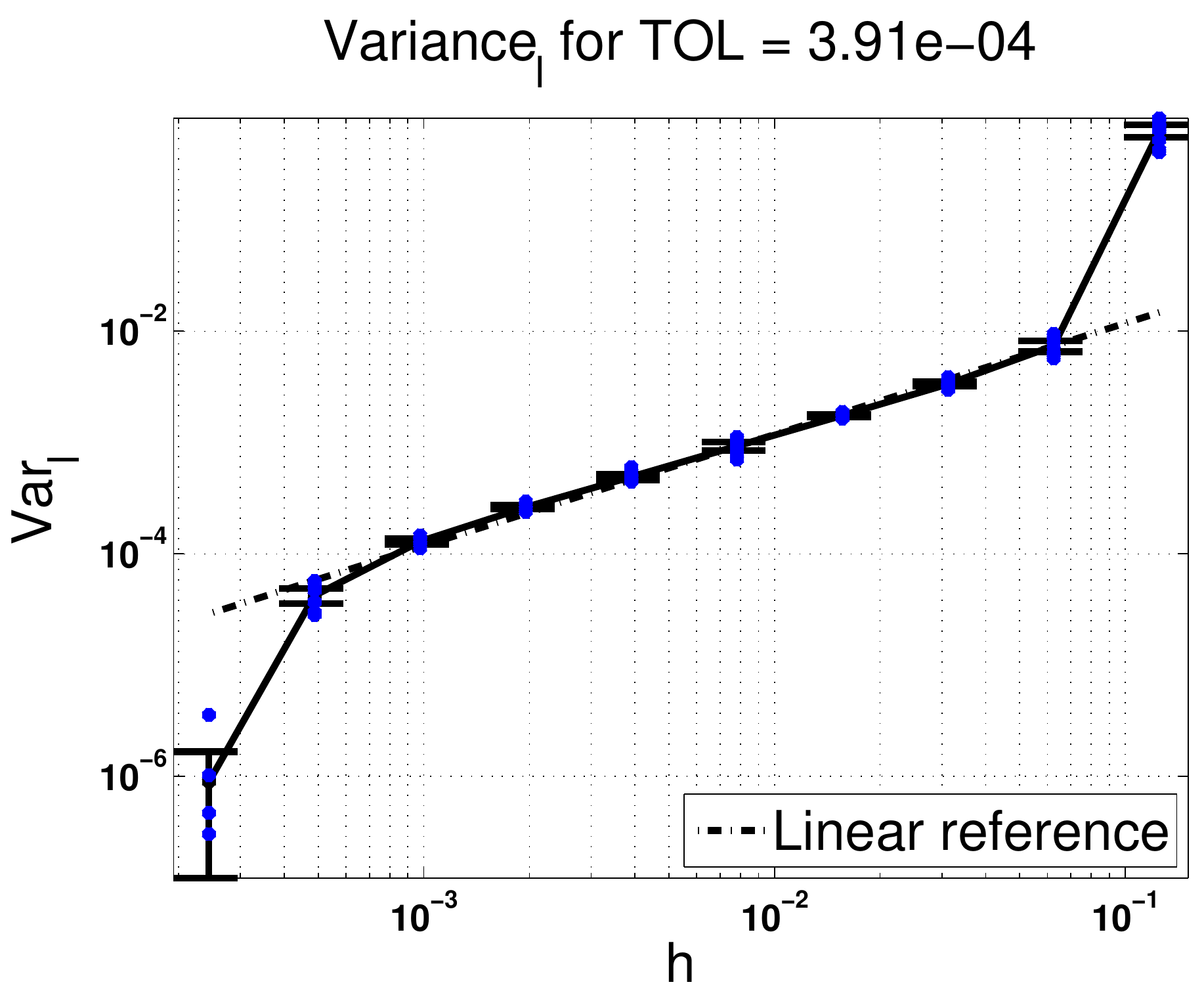}
\end{minipage}
\begin{minipage}{0.49\textwidth}
\includegraphics[scale=0.31]{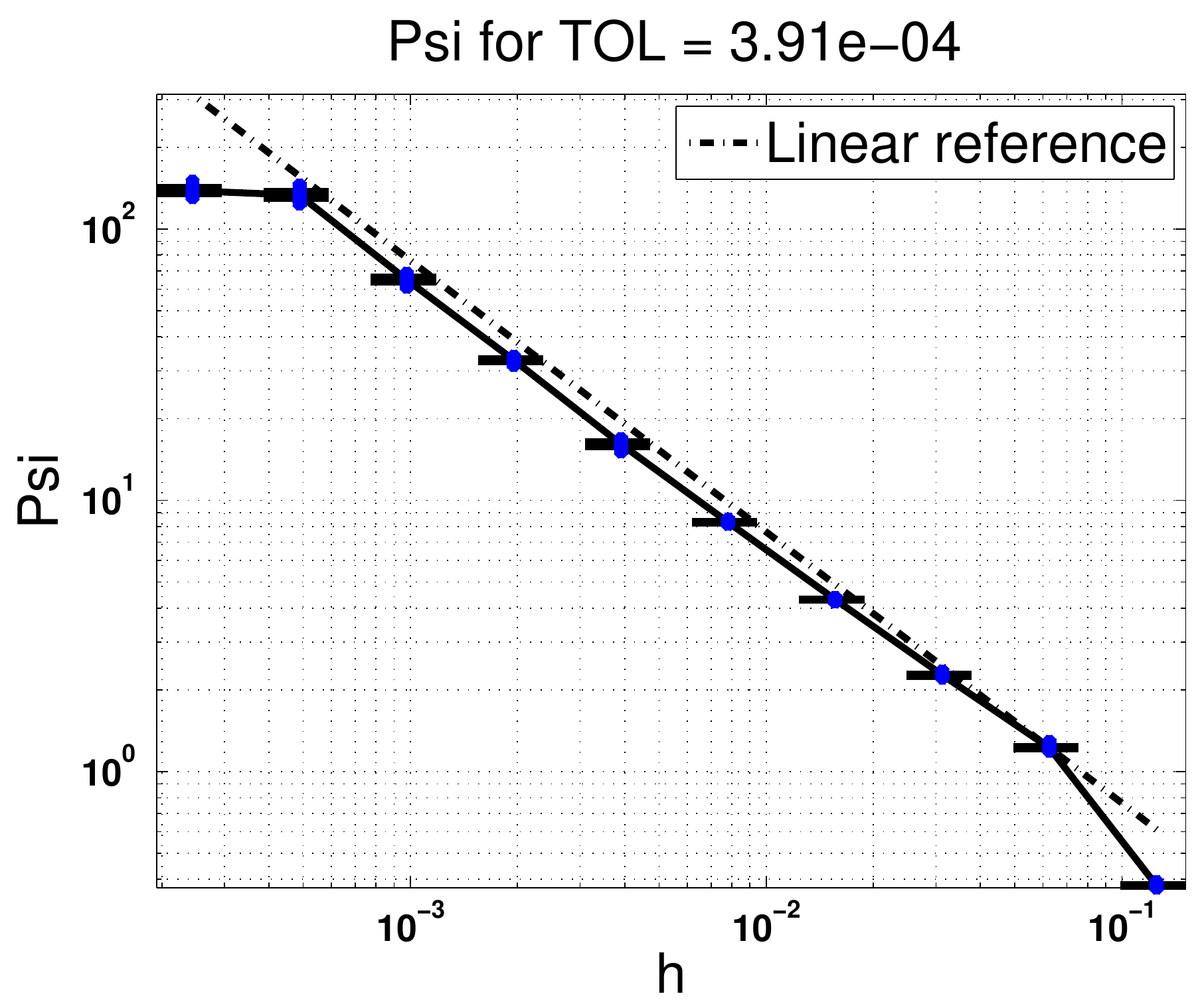}
\end{minipage}
\begin{minipage}{0.49\textwidth}
\includegraphics[scale=0.31]{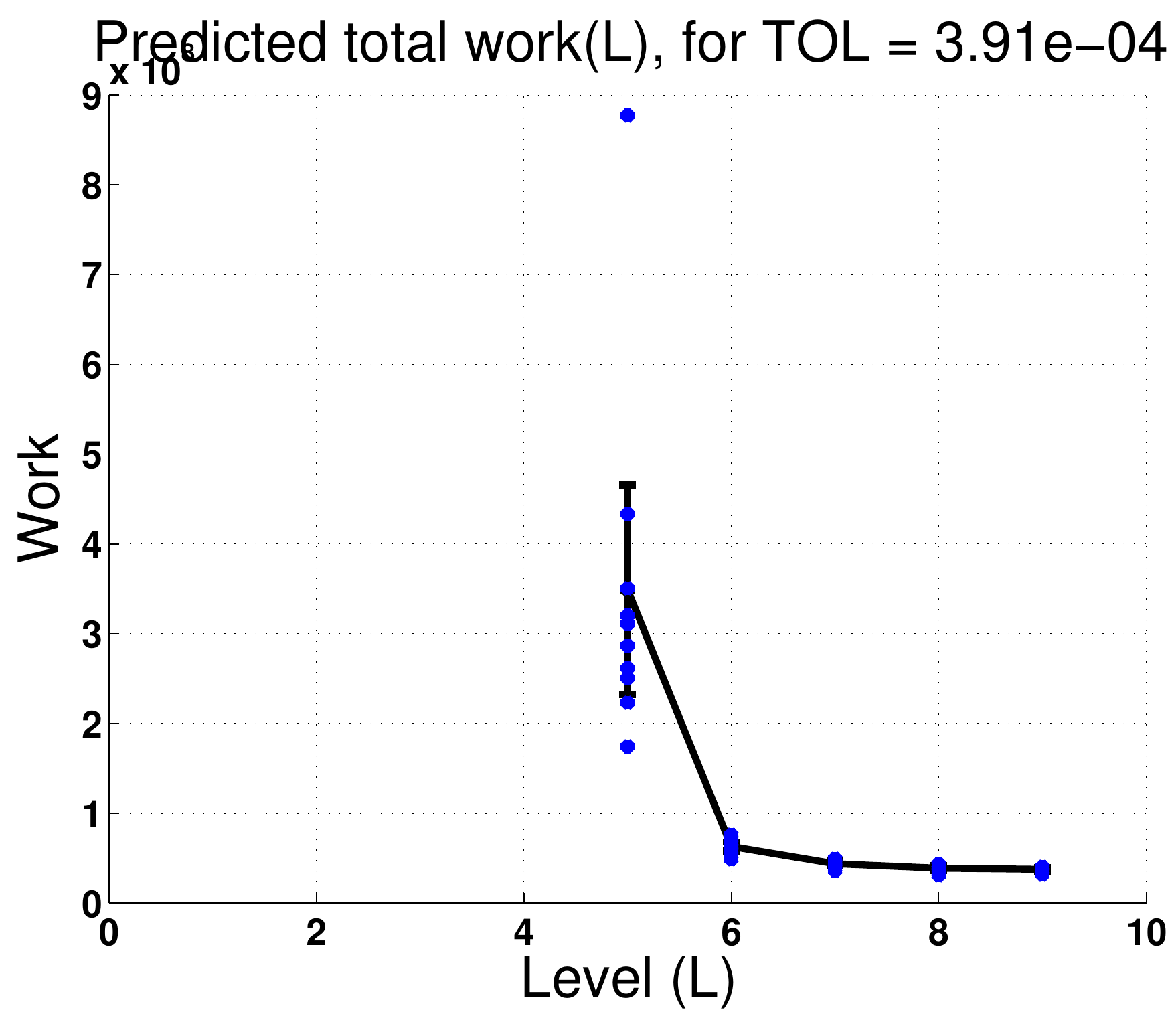}
\end{minipage}
\caption{Upper left: estimated weak error, $\WEH{\ell}$, as a function of the time mesh size, $h$. Upper right: estimated variance of the difference between two consecutive levels, $\hmV{\ell}$, as a function of $h$. Lower left: estimated path work, $\hat \psi_\ell$, as a function of $h$. Lower right: estimated total computational work, $\sum_{l=0}^L \hat \psi_l M_l$, as a function of the level, $L$.
}
\label{fig:dec2-diag}
\end{figure}

In Figure \ref{fig:dec2-diag}, we can observe how the estimated weak error, $\WEH{\ell}$, and the estimated variance of the difference of the functional between two consecutive levels, $\hmV{\ell}$, decrease linearly as we refine the time mesh, which corresponds to a tau-leap dominated regime. This linear relationship for the variance starts at level $1$, as expected. When the MNRM dominated regime is reached, both quickly converge to zero as expected.
The estimated total path work, $\hat \psi_{\ell}$, increases as we refine the time mesh. Observe that it increases linearly for the coarser grids, until it reaches a plateau, which corresponds to the pure MNRM case where the computational cost is independent of the grid size. 
In the lower right panel, we show the total computational work, only for the cases in which $\WEH{\ell} < TOL{-}TOL^2$. 
\newcommand{\TOLDEC}{3.91e-04}

\begin{figure}[h!]
\centering
\begin{minipage}{0.49\textwidth}
\includegraphics[scale=0.31]{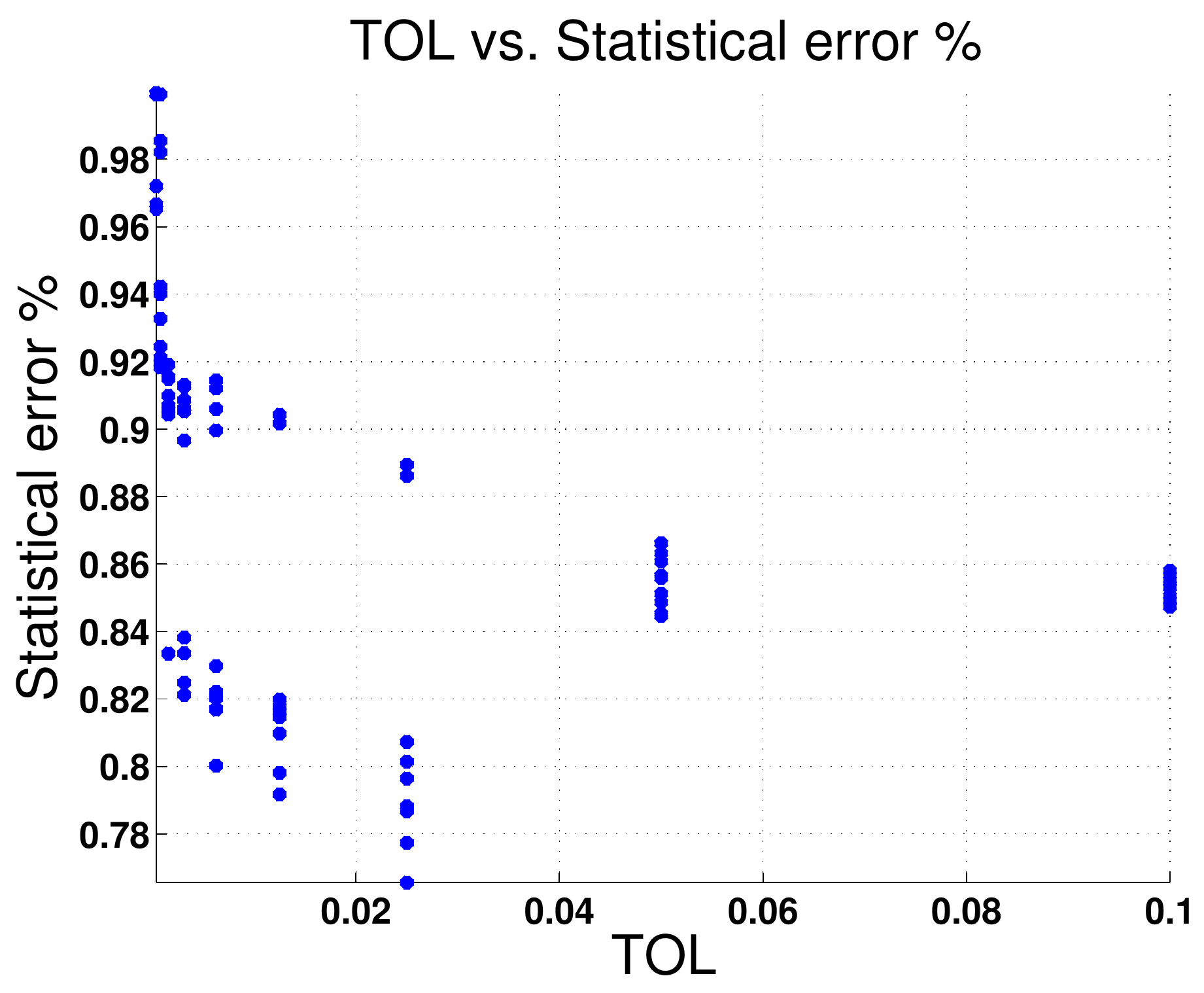}
\end{minipage}
\begin{minipage}{0.49\textwidth}
\includegraphics[scale=0.31]{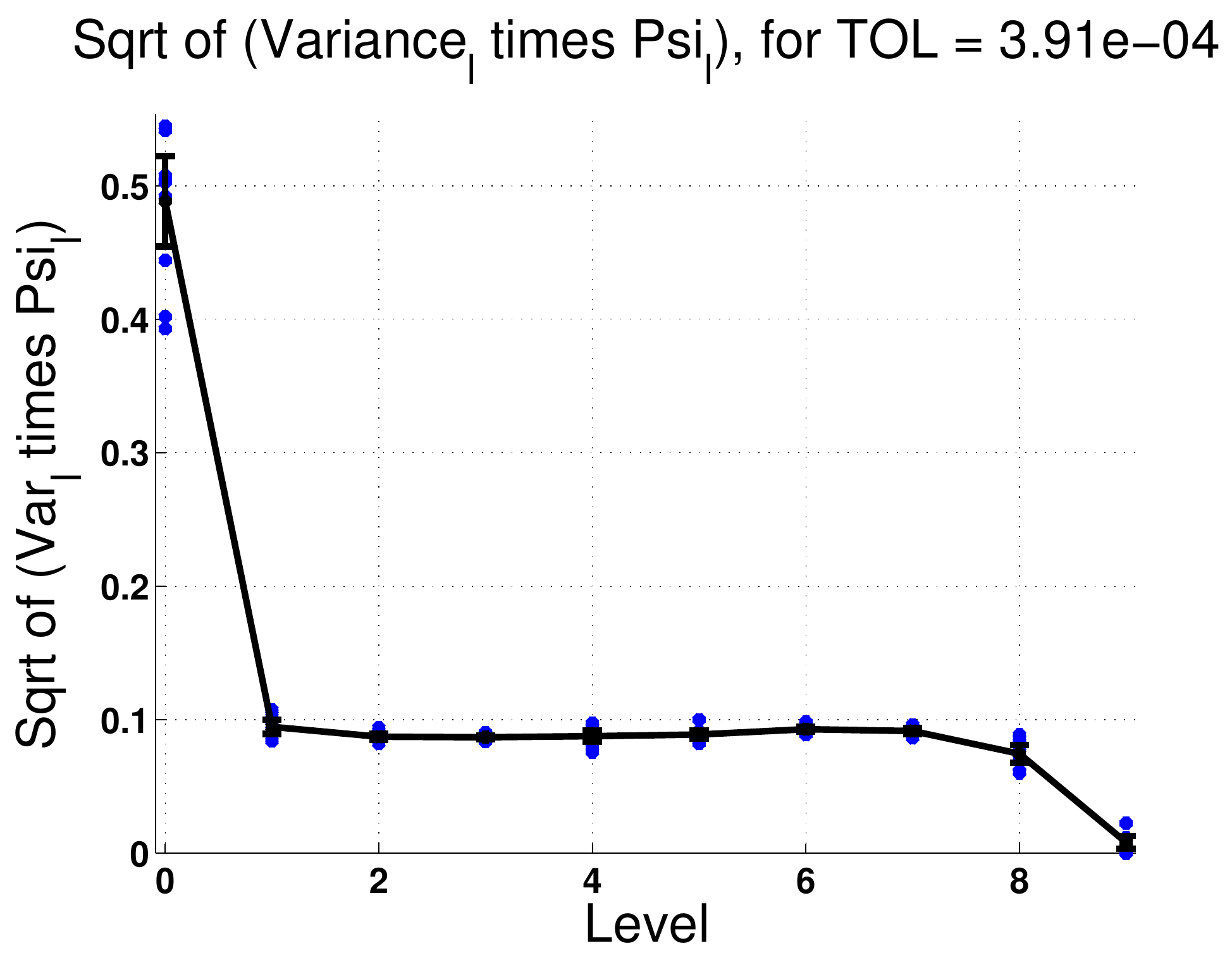}
\end{minipage}
\caption{Left: Percentage of the statistical error over the total error. As we mentioned in Section \ref{sec:complex}, it is well above $0.5$ for all the tolerances. Right: $\sqrt{\hmV{\ell}\hat{\psi}_\ell}$, as a function of $\ell$, for the smallest tolerance, which  decreases as the level increases. Observe that the contribution of level 0 is less than 50\% of the sum of the other levels.
}
\label{fig:statdec2}
\end{figure}

\begin{figure}[h!]
\centering
\begin{minipage}{0.49\textwidth}
\includegraphics[scale=0.31]{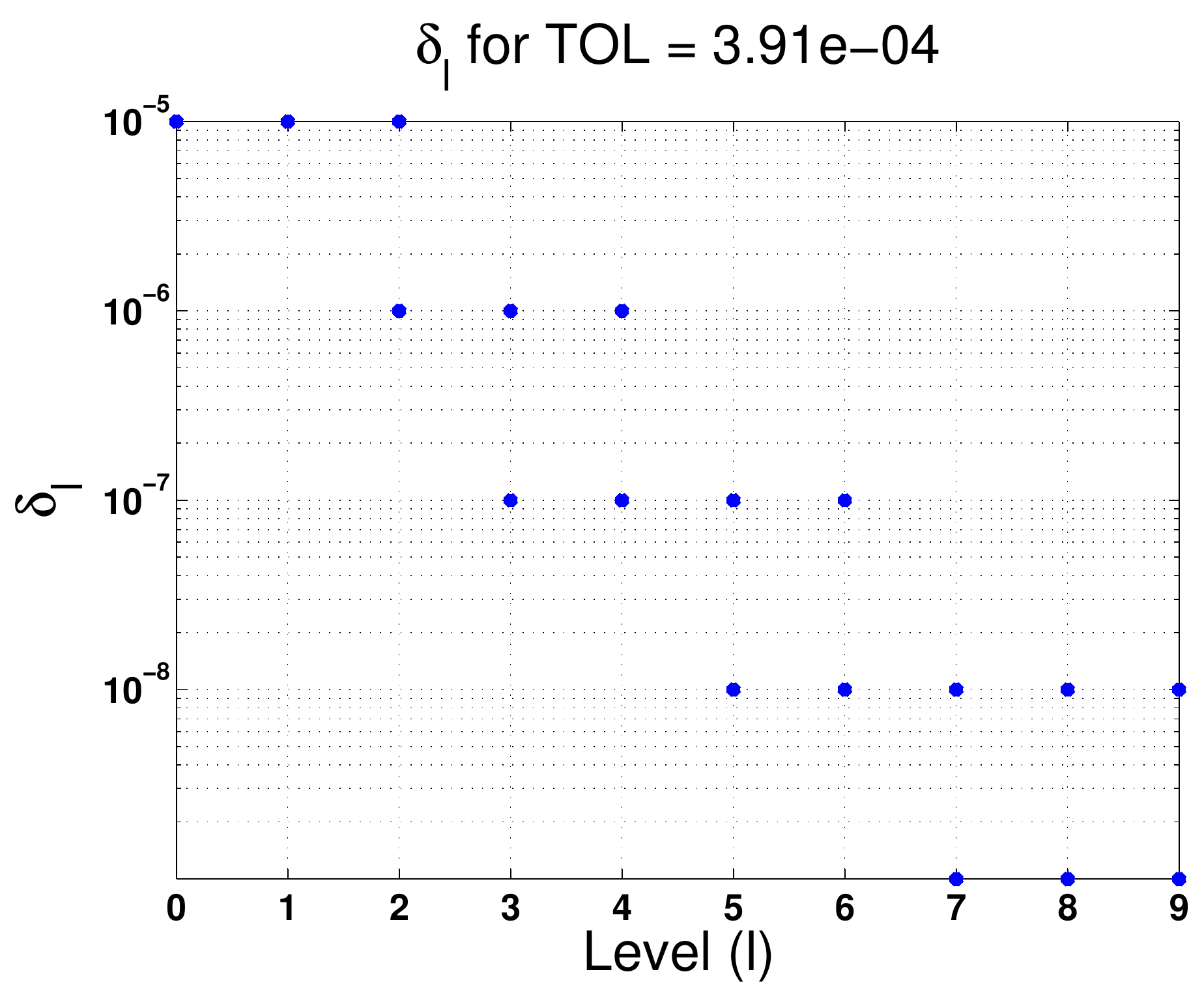}
\end{minipage}
\begin{minipage}{0.49\textwidth}
\includegraphics[scale=0.31]{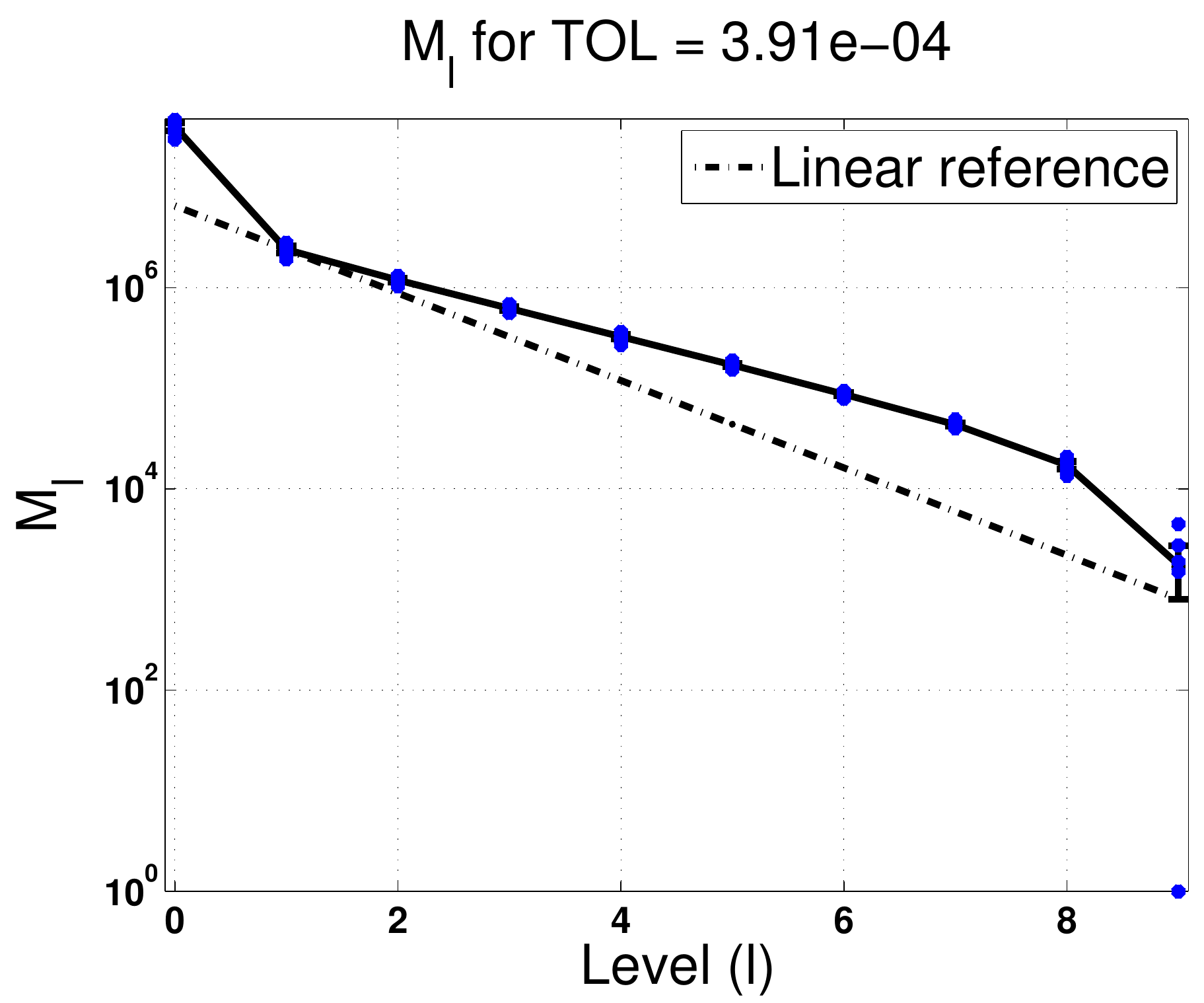}
\end{minipage}
\caption{The one-step exit probability bound, $\delta_\ell$, and $M_\ell$ for $\ell{=}0,1,...,L^*$, for the smallest tolerance.
}
\label{fig:dec2-out}
\end{figure}

In Figure \ref{fig:dec2-out}, we show the main outputs of the phase II algorithm, $\delta_\ell$ and $M_\ell$ for $\ell=0,...,L^*$, for the smallest considered tolerance. In this example, $L^*$ is 8 or 9, depending on the run. We observe that the number of realizations decreases slower than linearly, from levels $1$ to $L^*{-}1$, until it drops, due to the change to a MNRM dominated regime. 

\begin{figure}[h!]
\centering
\begin{minipage}{0.49\textwidth}
\includegraphics[scale=0.31]{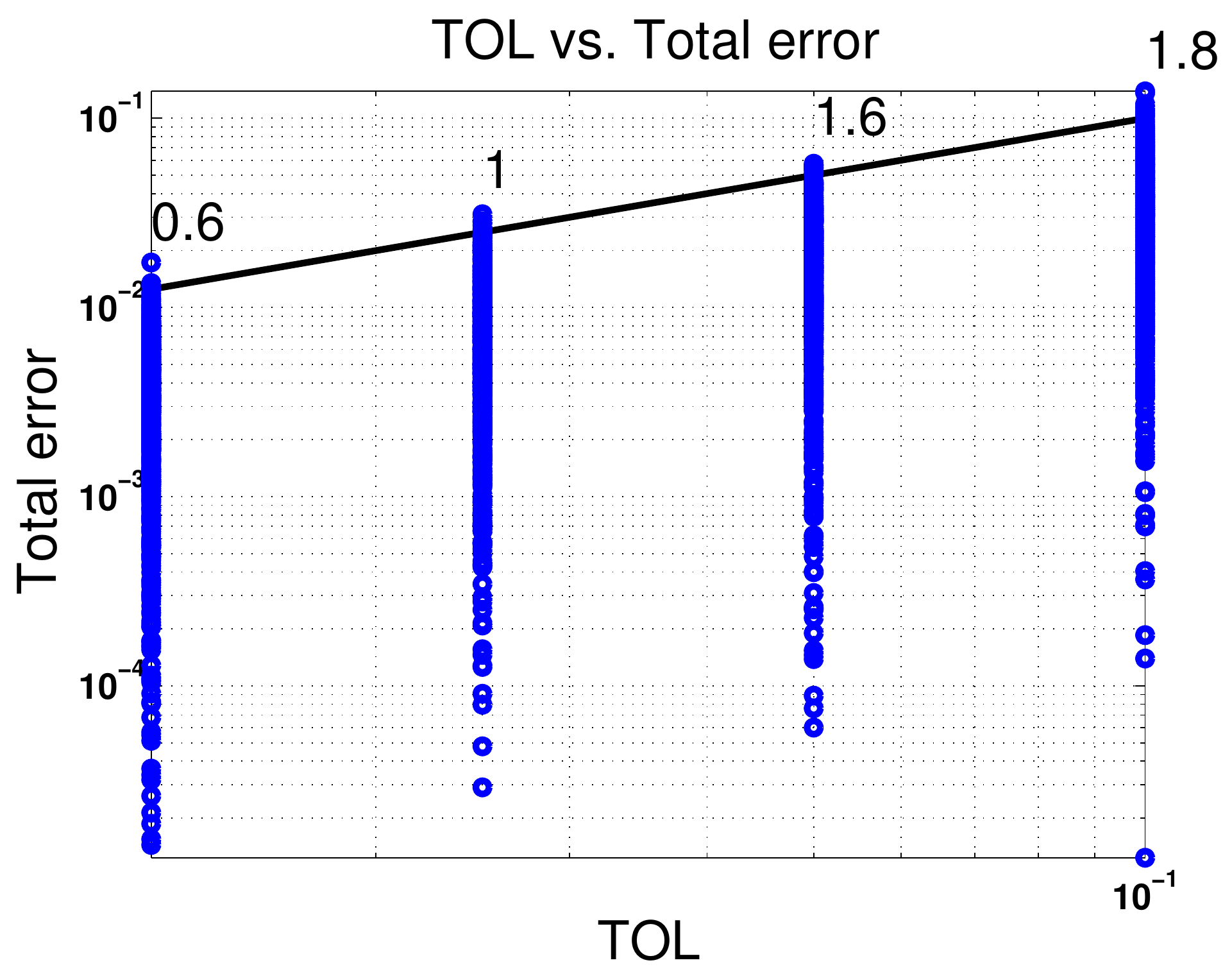}
\end{minipage}
\begin{minipage}{0.49\textwidth}
\includegraphics[scale=0.31]{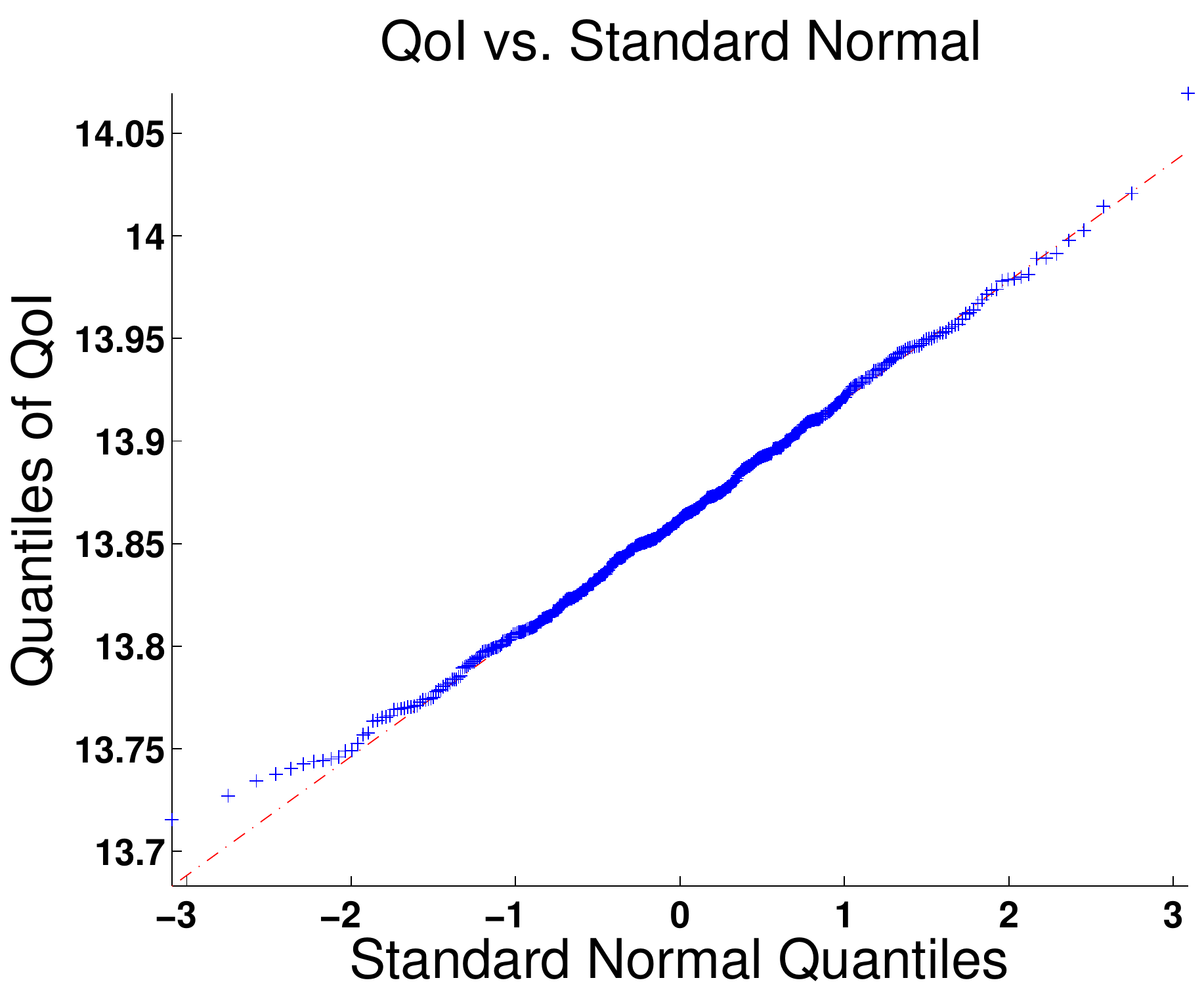}
\end{minipage}
\caption{Left: $TOL$ versus the actual computational error. The numbers above the straight line show the percentage of runs that had errors larger than the required tolerance. We observe that in all cases, the computational error follows the imposed tolerance with the expected confidence of 95\%. Right: quantile-quantile plot based on realizations of $\mathcal{M}_L$.
}
\label{fig:effvir}
\end{figure}

In Figure \ref{fig:effvir}, 
we show $TOL$ versus the actual computational error. It can be seen that the prescribed tolerance is achieved with the required confidence of 95\%, since $C_A{=}1.96$, for all the tolerances.
The QQ-plot in the right part of Figure \ref{fig:effvir} was obtained as follows: i) for the range of tolerances specified in the first column of Table \ref{tab:vir}, we ran the phase II algorithm 5 times, ii) for each output of the calibration algorithm, we sampled the multilevel estimator $\mathcal{M}_L$, defined in \ref{MLMCest}, 100 times. This plot reaffirms our assumption about the Gaussian distribution of the statistical error. 

\begin{rem}
In the simulations, we observe that, as we refine $TOL$, the optimal number of levels approximately increases logarithmically, which is a desirable feature. We fit the model $L^*=a\log(TOL^{-1})+b$, obtaining $a{=}1.47$ and $b{=}3.56$.
\end{rem}

\begin{rem}[Pareto rule]
\label{rem:pareto}
Using the cost-based rule (see remark \ref{rem:pare}), we estimate the threshold for the Pareto rule, obtaining $\nu = 0.95419$. It turns out that, for this example, $\hat W_{MixPareto}/ \hat W_{Mix}$ ranges from $0.6$ to $0.75$ (for most $TOL$s). This shows that it is possible to  increase the computational work gains further in some examples. 
\end{rem}

\begin{rem}
The savings in computational work when generating Poisson random variables heavily depend on MATLAB's performance capabilities. In fact, we would expect better results from our method if we were to implement our algorithms in more performance-oriented languages or if we were to sample Poisson random variables in batches.
\end{rem}

\subsection*{A Simple Stiff System}
This model, adapted from \cite{gillespie2005}, has three species and a mixture of fast and slow reaction channels,
$$X_1 \xrightleftharpoons[c_2]{c_1} X_2 \xrightarrow{c_3} X_3 \xrightarrow{c_4} \emptyset, \,\,\, c_2 \gg c_3 > c_4 \PERIOD$$

\begin{figure}[h!]
\centering
\begin{minipage}{0.49\textwidth}
\includegraphics[scale=0.34]{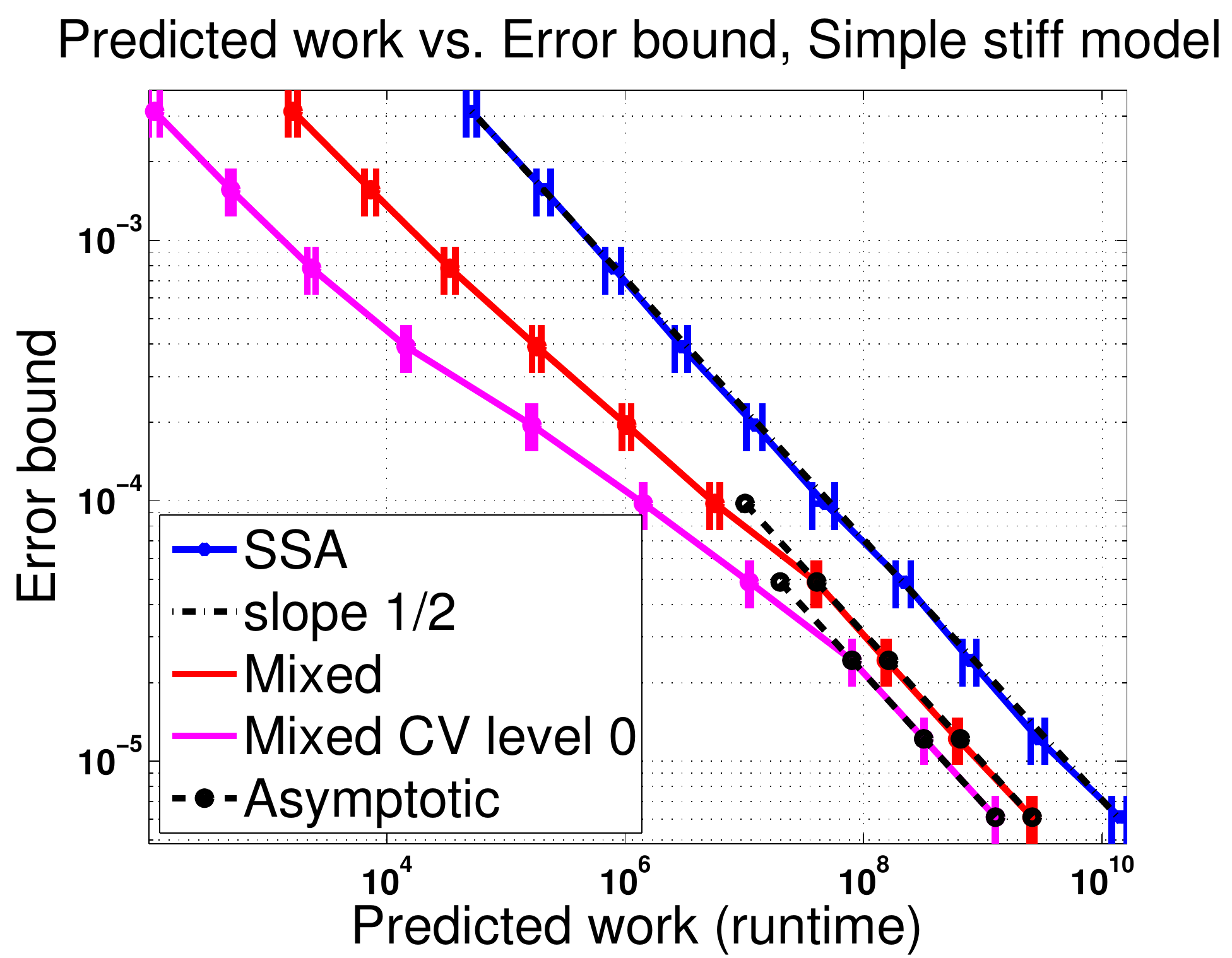}
\end{minipage}
\hfill
\begin{minipage}{0.49\textwidth}
\begin{small}
\begin{tabular}{lccc}
$TOL$ & $L^*$ & $\frac{\hat{W}_{ML_{cv}}}{\hat{W}_{\ssa}}$ & $\frac{\hat{W}_{ML}}{\hat{W}_{\ssa}}$ \\ \noalign{\smallskip} \hline\noalign{\smallskip} 
3.13e-03 & 1.0  & 0.002 $\pm$0.0004 & 0.03 $\pm$0.001\\ 
1.56e-03 & 1.0  & 0.003 $\pm$0.0004 & 0.04 $\pm$0.001 \\ 
7.81e-04 & 1.0  & 0.003 $\pm$0.0010 & 0.04 $\pm$0.002\\ 
3.91e-04 & 1.0  & 0.004 $\pm$0.0004 & 0.06 $\pm$0.003\\ 
1.95e-04 & 2.0  & 0.013 $\pm$0.0015 & 0.09 $\pm$0.008\\ 
9.77e-05 & 3.0  & 0.027 $\pm$0.0040 & 0.13 $\pm$0.016\\ 
4.88e-05 & 4.0  & 0.065 $\pm$0.0146 & 0.19 $\pm$0.025 \\ 
2.44e-05 & 6.0  & 0.100 $\pm$0.0136 & 0.21 $\pm$0.020\\ 
1.22e-05 & 6.0  & 0.109 $\pm$0.0299 & 0.22 $\pm$0.029\\ 
6.10e-06 & 6.0  & 0.108 $\pm$0.0168 & 0.19 $\pm$0.020\\ 
\noalign{\smallskip}\hline 
\end{tabular}
\end{small}
\end{minipage}
\caption{Left: Predicted work (runtime) versus the estimated error bound, with $95\%$ confidence intervals, for the simple stiff model with and without using the control variate at level 0, as described in Section \ref{sec:cvar}. 
Right: Details of the ensemble run of the phase II algorithm using the control variate (third column) and without using the control variate (fourth column).
As an example, the fifth row of the table tells us that, for a tolerance $TOL{=}1.95\cdot 10^{-4}$, 2 levels are needed on average. 
The work of the multilevel mixed method using the control variate at level 0 is, on average, $1\%$ of the work of the SSA. When not using the control variate, it is $9\%$.
Confidence intervals at 95\% are also provided.
}
\label{fig:stiff-worktime}
\end{figure}

Its stoichiometric matrix and propensity functions, $a_j:\zset_+ \rightarrow \rset$, are given by 
\begin{align*}
\nu =\left( 
 \begin{array}{ccc}     -1  &   1 &    0 \\
     1  & -1 &   0\\
     0 &   -1  &   1\\
     0  &   0  &  -1
 \end{array} 
 \right) ^{tr}\mbox{    and     }  a(X) =\left( \begin{array}{c} c_1 X_1 \\ 
 c_2 X_2 \\ c_3 X_2 \\ c_4 X_3 \end{array} \right)\COMMA
\end{align*}
respectively, where $g(X(t))=X_3(t)$.
In this model, successive firings of the reaction $X_2 \rightarrow X_3$ are separated by many reversible firings between $X_1$ and $X_2$, which takes a lot of computational work in a standard SSA run. In \cite{gillespie2013}, Gillespie et al. claim that this inefficiency cannot be addressed using ordinary tau-leaping because of the stiffness of the system. We show here that we have substantial gains using our mixed method, which also controls the global error. In this example, we also show the performance of the control variate idea, presented in Section \ref{sec:cvar}.
We analyze 10 independent runs of the phase II algorithm (see Section \ref{sec:estim}), using different relative tolerances. 
In Figure \ref{fig:stiff-worktime}, we show the total predicted work (runtime) for the multilevel mixed method with and without a control variate at level 0 and for the SSA method versus the estimated error bound. We also show the estimated asymptotic work of the multilevel mixed method. 
Observe that, for practical tolerances, the computational work gains with respect to the SSA method, when using the control variate, are of a factor of 500 times. Without using the control variate, computational gains are also substantial. 

\section{Conclusions}
\label{sec:conclusions}

In this work, we addressed the problem of approximating the quantity of interest $\expt{g(X(T))}$, where $X$ is a non-homogeneous Poisson process that  describes a stochastic reactions network, and $g$ is a given suitable observable of $X$, within a given prescribed relative tolerance, $TOL{>}0$, up to a given confidence level at near-optimal computational work.

We developed an automatic, adaptive reaction-splitting multilevel Monte Carlo method, based on our Chernoff tau-leap mehthod \cite{ourSL,ourML}. 
Its computational complexity is $\Ordo{TOL^{-2}}$. This method can be therefore seen as a variance reduction of the SSA, which has the same complexity. 
In our numerical examples, we obtained substantial gains with respect to SSA and, for systems in which the set of reaction channels can be adaptively partitioned into ``high'' and ``low'' activity, over our previous multilevel hybrid Chernoff tau-leap method.

We also presented a novel control variate for $g(X(T))$, which adds negligible computational cost when simulating a path of $X(T)$, and it may lead to additional dramatic cost reductions.

\section*{Acknowledgments}
The research reported here was supported by King Abdullah University of Science and Technology (KAUST).
The authors are members of the SRI Center for
Uncertainty Quantification in Computational Science and Engineering at KAUST. 
\appendix



\end{document}